\documentclass[a4paper,11pt,twoside]{article}
\usepackage{pst-fill,pst-grad}
\usepackage{textcomp}
\usepackage[english]{babel}
\usepackage{graphicx}
\usepackage{amsmath}
\usepackage{float}
\usepackage[matrix,arrow,curve]{xy}
\usepackage{pstricks}
\usepackage{mathtools}
\usepackage{amsmath,amsfonts,verbatim,afterpage,theorem,euscript,mathrsfs,amssymb}
\usepackage{amsfonts}
\usepackage{amssymb}
\usepackage{array}
\usepackage{dsfont}
\usepackage{hyperref}
\def \vu{\vec{u}}
\def \vU{\vec{U}}
\def \vpsi{\vec{\psi}}
\def \vf{\vec{f}}
\def \div{\mathrm{div}}
\def \vn{\vec{\nabla}}

\newtheorem{Theoreme}{Theorem}

\newtheorem{Proposition}{Proposition}[section]
\newtheorem{Lemme}{Lemma}[section]

\newtheorem{Remarque}{Remark}[section]

\numberwithin{equation}{section}
\title{\bf Global mild solutions in a critical setting for a forced fractional Boussinesq system}

\usepackage{authblk}
\author{Diego Chamorro\footnote{\emph{diego.chamorro@univ-evry.fr} }}
\affil{\footnotesize LaMME, Univ. Evry, CNRS, Universit\'e Paris-Saclay, 91025, Evry, France.}
\author{Maxence Mansais\footnote{\emph{maxence.mansais@ens-paris-saclay.fr} }}
\affil{\footnotesize LaMME, Univ. Evry, CNRS, Universit\'e Paris-Saclay, 91025, Evry, France.}
\begin{document}
\maketitle
\begin{abstract}
\noindent We study here mild solutions for the forced, incompressible fractional Boussinesq system. Under suitable estimates for the terms involved (in an adapted functional framework) we can invoque a fixed point argument in order to obtain mild solutions. Although many functional spaces can be considered, we are interested here in a critical setting which ensures the existence of global solutions and we will work in particular with parabolic Morrey spaces which provide one of the largest critical functional frameworks available for constructing mild solutions for the fractional Boussinesq equations.
\end{abstract}

\section{Introduction and presentation of the results}
In this article we are interested in studying mild solutions for the fractional incompressible Boussinesq system over $\mathbb{R}^d$ with $d\geq 2$, which is given by the following coupled equations
\begin{equation}\label{EquationNSB_intro}
\begin{cases}
\partial_t \vec{u}  = -(-\Delta)^\frac{\alpha}{2} \vec{u}-\mathrm{div} \left(\vec{u}\otimes\vec{u}\right) -\vn p + \theta \vec{e}_d + \vec{f}, \qquad \mathrm{div}(\vec{u})=0, \qquad (1<\alpha< 2),\\[3mm]
\partial_t\theta = -(-\Delta)^\frac{\alpha}{2}\theta + \mathrm{div}(\theta \ \vec{u})+g,\\[3mm]
\vec{u}(0,\cdot)  = \vec{u}_0,\quad \mathrm{div}(\vec{u}_0)=0, \qquad \theta(0,\cdot)=\theta_0,
\end{cases}
\end{equation}
where $\vec{u}_0:\mathbb{R}^d\longrightarrow \mathbb{R}^d$ denotes the initial velocity of the fluid,  $\theta_0:\mathbb{R}^d\longrightarrow \mathbb{R}$ is the initial temperature of the fluid, $\vec{u}:[0,+\infty[\times \mathbb{R}^d\longrightarrow \mathbb{R}^d$ is the velocity field which is assumed to be divergence free (\emph{i.e.} $\mathrm{div}(\vec{u})=0$), $p:[0,+\infty[\times \mathbb{R}^d\longrightarrow \mathbb{R}$ is the internal pressure and $\theta:[0,+\infty[\times \mathbb{R}^d\longrightarrow \mathbb{R}$ is the temperature of the fluid (which can be viewed as convection force which is both diffused and transported along the flow of the velocity $\vu$). The functions $\vec{f}:[0,+\infty[\times \mathbb{R}^d\longrightarrow \mathbb{R}^d$ and $g:[0,+\infty[\times \mathbb{R}^d\longrightarrow \mathbb{R}$ represent two given external forces. We recall that for an index $\alpha > 0$, the fractional Laplacian operator $(-\Delta)^\frac{\alpha}{2}$ can be defined as a Fourier multiplier through the expression $\widehat{(-\Delta)^\frac{\alpha}{2}\phi}(\xi)= |\xi|^\alpha \widehat{\phi}(\xi)$, where $\phi:\mathbb{R}^d\longrightarrow \mathbb{R}$ is a suitable function (say $\phi\in \mathcal{S}(\mathbb{R}^d)$). For vectorial functions $\vec{\phi}:\mathbb{R}^d\longrightarrow \mathbb{R}^d$ the expression $(-\Delta)^\frac{\alpha}{2}\vec{\phi}$ is considered component-wise.\\

The classical Boussinesq system (\emph{i.e.} when $\alpha=2$) has been studied by multiple authors. For weak solutions \emph{à la Leray} and their properties see \cite{BrandoleseSchonbek}, \cite{FanZhou}, \cite{GuoYuan} and for mild solutions \emph{à la Kato} see \cite{BrandoleseHe}, \cite{DanchinPaicu}, \cite{CannonDiBenedetto} as well as the references therein. For other results related to this system (in 3D or in 2D) see \cite{Chae}, \cite{GancedoGarcia}, \cite{HmidiZerguine}, \cite{Sawada}, \cite{YeXuXue}. The fractional Boussinesq system (\ref{EquationNSB_intro}) has been a little less studied than the classical equations, see however the articles \cite{JiuMiao}, \cite{HmidiKR}, \cite{Stefanov} and \cite{YangJ} for some results in the fractional framework.\\ 

In this article we are particularly interested in constructing mild solutions of the system (\ref{EquationNSB_intro}) in a very large critical functional setting (essentially based on parabolic Morrey spaces and this approach seems to be new for the fractional Boussinesq equations) and in order to present our framework we need to recall some basic facts.\\

Indeed, we first remark that it is possible to link the study of the pressure $p$ to the velocity field $\vu$ and to the temperature $\theta$: by formally applying the divergence to the first equation of the system (\ref{EquationNSB_intro}) and using the divergence free property of the velocity field ($\div(\vu)=0$) we obtain the following expression
$$(-\Delta) p=\mathrm{div}\left(\mathrm{div}(\vec{u}\otimes \vec{u}) - \vec{f} - \theta\vec{e}_d\right),$$
which allows us to deduce information over the pressure as long as we have suitable informations over the velocity field $\vu$, the temperature $\theta$ and the external force $\vf$. This enables us to consider the functions $\vu$ and $\theta$ as the main variables of the system (\ref{EquationNSB_intro}) and we will thus focus our study in these functions. To continue, and in order to get rid of the pressure, we consider now the Leray projector defined by the expression $\mathbb{P}(\vpsi)=\vpsi + \vn (-\Delta)^{-1}\div(\vpsi)$ and since we have $\mathbb{P}(\vn p)=0$ as well as $\mathbb{P}(\vu)=\vu$ (since $\div(\vu)=0$), if we apply this projector to the first equation of (\ref{EquationNSB_intro}), we obtain the equation
$$\partial_t \vec{u}  = -(-\Delta)^\frac{\alpha}{2} \vec{u}-\mathbb{P}\left(\mathrm{div} \left(\vec{u}\otimes\vec{u}\right)\right) + \mathbb{P}(\theta\vec{e}_d) + \mathbb{P}(\vec{f}),$$
where the pressure $p$ does not intervene. We will thus work with the following system 
\begin{equation}\label{EquationNSB_SansPression}
\begin{cases}
\partial_t \vec{u}  = -(-\Delta)^\frac{\alpha}{2} \vec{u}-\mathbb{P}\left(\mathrm{div} \left(\vec{u}\otimes\vec{u}\right)\right) + \mathbb{P}(\theta\vec{e}_d) + \mathbb{P}(\vec{f}),\\[3mm]
\partial_t\theta = -(-\Delta)^\frac{\alpha}{2}\theta + \mathrm{div}(\theta \ \vec{u})+g,
\end{cases}
\end{equation}
and a very important feature of this system is the following: if $(\vu, \theta)$ is a solution of (\ref{EquationNSB_SansPression}), then for all $\lambda>0$, the rescaled couple $(\vu_\lambda, \theta_\lambda)$ where 
\begin{equation}\label{Def_Scaling}
\vu_\lambda(t,x)=\lambda^{\alpha-1}\vu(\lambda^\alpha t, \lambda x),\quad \mbox{  }\quad \theta_\lambda(t,x)=\lambda^{2\alpha-1}\theta(\lambda^\alpha t ,\lambda x),
\end{equation}
is also a solution of the same system, of course with the rescaled forces 
\begin{equation}\label{Def_Scaling_Forces}
\vf_\lambda(t,x)=\lambda^{2\alpha-1}\vf(\lambda^\alpha t ,\lambda x),\quad \mbox{  }\quad g_\lambda(t,x)=\lambda^{3\alpha-1}g(\lambda^\alpha t ,\lambda x),
\end{equation}
and with the rescaled initial data:
\begin{equation}\label{Def_Scaling_DonneeInitiales}
\vu_{\lambda, 0}(x)= \lambda^{\alpha-1}\vu_0(\lambda x),\quad \mbox{  }\quad \theta_{\lambda, 0}(x)=\lambda^{2\alpha-1}\theta_0(\lambda x).
\end{equation}
We soon shall see how this homogeneity (or scaling) property of the system (\ref{EquationNSB_SansPression}) will be essential in what follows. \\

As we are mainly interested in studying mild solutions for the system (\ref{EquationNSB_SansPression}) (see Theorem \ref{Theorem_Triebel_Morrey} below), we will consider its integral formulation:
\begin{equation}\label{Formulation_Integrale}
\begin{split}
&\vec{u}=\mathfrak{p}_t\ast \vec{u_0}-\int_0^t \mathfrak{p}_{t-s}\ast\mathbb{P}\left(\mathrm{div}(\vec{u}\otimes \vec{u})\right)ds + \int_0^t \mathfrak{p}_{t-s}\ast\mathbb{P}(\theta \vec{e}_d)ds + \int_0^t \mathfrak{p}_{t-s}\ast\mathbb{P}(\vec{f})ds\\
&\theta =\mathfrak{p}_t\ast \theta_0+\int_0^t \mathfrak{p}_{t-s}\ast \mathrm{div}(\theta \vec{u})ds+ \displaystyle{\int_0^t} \mathfrak{p}_{t-s}\ast g\,ds,
\end{split}
\end{equation}
where $\mathfrak{p}_t$, with $t>0$ is the fractional heat kernel associated to the semi-group $e^{-t(-\Delta)^{\frac{\alpha}{2}}}$ (we thus have in the Fourier level the expression $\widehat{\mathfrak{p}_t}(\xi)=e^{-t|\xi|^\alpha}$).\\


Before we state our first result, we need to introduce some definitions.
\begin{itemize}
\item {\bf Morrey spaces}. We will say that a locally integrable function $\phi:\mathbb{R}^d\longrightarrow \mathbb{R}$ belongs to the Morrey space $M^{p,q}(\mathbb{R}^d)$ with $1\leq p\leq q<+\infty$, if we have 
\begin{equation}\label{Def_EspaceMorrey}
\|\phi\|_{M^{p,q}}=\underset{x\in \mathbb{R}^d}{\sup}\,\underset{r>0}{\sup}\;\frac{1}{r^{d(\frac{1}{p}-\frac{1}{q})}}\left(\int_{B(x,r)}|\phi(y)|^pdy\right)^{\frac{1}{p}}<+\infty.
\end{equation}
Morrey spaces are a generalization of Lebesgue spaces and from this expression we easily obtain the identification $M^{p,p}(\mathbb{R}^d)=L^p(\mathbb{R}^d)$. See more details on these functional spaces in the books \cite{Rafeiro} and \cite{SickelYangYuan}.
\item {\bf Triebel-Lizorkin-Morrey spaces}. We will need the following thermic characterization of Triebel-Lizorkin spaces based on Morrey spaces. Indeed, for a regularity index $\sigma <0$ and for two parameters $1\leq p\leq q<+\infty$, we will say that a function $\phi:\mathbb{R}^d\longrightarrow \mathbb{R}$ belongs to the Triebel-Morrey spaces $T_{M^{p,q}}^\sigma(\mathbb{R}^d)$ if 
\begin{equation}\label{Def_EspaceTriebelLizorkinMorrey}
\|\phi\|_{T_{M^{p,q}}^\sigma}=\left\|\left(\int_{0}^{+\infty}t^\frac{-\sigma p}{\alpha}|\mathfrak{p}_t\ast \phi(\cdot)|^p\frac{dt}{t}\right)^{\frac{1}{p}}\right\|_{M^{p,q}}<+\infty,
\end{equation}
where $\mathfrak{p}_t$ stands for the kernel associated to the semi-group $e^{-t(-\Delta)^{\frac{\alpha}{2}}}$ with $1<\alpha<2$ and where $M^{p,q}(\mathbb{R}^d)$ is a Morrey space given in the expression (\ref{Def_EspaceMorrey}) above.
\item {\bf Parabolic Morrey spaces}. Now, for a function $\psi:[0,+\infty[\times\mathbb{R}^d\longrightarrow \mathbb{R}$ we will say that it belongs to the homogeneous \emph{parabolic} Morrey space $\mathcal{M}^{p,q}_\alpha([0,+\infty[\times\mathbb{R}^d)$ where $1\leq p\leq q<+\infty$ and $1<\alpha<2$ if 
\begin{equation}\label{Definition_MorreyPara}
\|\psi\|_{\mathcal{M}_\alpha^{p,q}}=\underset{r>0}{\mathrm{sup}} \ \underset{(t,x)\in [0,+\infty[ \times \mathbb{R}^d}{\mathrm{sup}} \ \frac{1}{r^{(d+\alpha)(\frac{1}{p}-\frac{1}{q})}} \left(\displaystyle{\iint_{\{|t-s|^\frac{1}{\alpha}+|x-y|<r\}}}|\psi(s,y)|^pdyds \right)^\frac{1}{p} <+\infty.
\end{equation}
It is easy to see that we have the space inclusion $\mathcal{M}^{p_1,q}_\alpha([0,+\infty[\times\mathbb{R}^d)\subset \mathcal{M}^{p_0,q}_\alpha([0,+\infty[\times\mathbb{R}^d)$ for all $1\leq p_0\leq p_1\leq q<+\infty$.
\item {\bf Parabolic Sobolev-Morrey spaces}. Finally, for a function $\psi:[0,+\infty[\times\mathbb{R}^d\longrightarrow \mathbb{R}$, we will say that it belongs to the homogeneous Sobolev-Morrey space $\dot{\mathcal{W}}^{-\gamma,p,q}_\alpha([0,+\infty[\times\mathbb{R}^d)$ where $1\leq p \leq q <+\infty$, $\gamma \in \mathbb{R}$ and $1< \alpha \leq 2$ if
\begin{equation}\label{Definition_SobolevMorrey}
\|\psi\|_{\dot{\mathcal{W}}^{-\gamma,p,q}_\alpha} = \|(-\Delta)^{-\frac{\gamma}{2}}\psi\|_{\mathcal{M}^{p,q}_\alpha} < +\infty. \\
\end{equation}
\end{itemize}
Of course these definitions can be generalized without any problem to vector valued functions. Note also that, in the definition of these parabolic spaces, the time interval $[0,+\infty[$ can be extended to $\mathbb{R}$ without any problem. See also Section \ref{Secc_Notations} below for more properties of these functional spaces. \\

With these functional spaces at hand we can now present our first result. 
\begin{Theoreme}[Initial data in Triebel-Lizorkin-Morrey spaces]\label{Theorem_Triebel_Morrey}
Let $d\geq2$ denote the dimension and consider a fractional power such that $1<\alpha < 2$. Fix now a real parameter $1\leq p<+\infty$ such that $4<\frac{3\alpha-2}{\alpha-1}<p\leq \frac{d+\alpha}{\alpha-1}$.
\begin{itemize}
\item[1)] If $\vu_0:\mathbb{R}^d\longrightarrow \mathbb{R}^d$ is a divergence free initial velocity such that $\vu_0\in T^{-\frac{\alpha}{p}}_{M^{p,q}}(\mathbb{R}^d)$ where we have $\frac{1}{q}=\frac{\alpha-1}{d+\alpha}-\frac{\alpha}{d}(\frac{1}{p}-\frac{\alpha-1}{d+\alpha})$,\\[2mm]  
if $\theta_0:\mathbb{R}^d\longrightarrow \mathbb{R}$ is an initial temperature such that $\theta_0\in T^{-\frac{\alpha}{\mathfrak{p}}}_{M^{\mathfrak{p}, \mathfrak{q}}}(\mathbb{R}^d)$ where $\mathfrak{p}=(\frac{\alpha-1}{2\alpha-1})p$ and $\frac{1}{\mathfrak{q}}=(\frac{2\alpha-1}{d+\alpha})-\frac{\alpha}{d}(\frac{1}{\mathfrak{p}}-\frac{2\alpha-1}{d+\alpha})$.\\

\item[2)] If $\vf:[0,+\infty[\times\mathbb{R}^d\longrightarrow \mathbb{R}^d$ is an external force such that 
$\vf \in \dot{\mathcal{W}}^{-\gamma,\mathfrak{m},\mathfrak{r}}_\alpha([0,+\infty[\times\mathbb{R}^d)$, where $0<\gamma <\alpha$ with $\mathfrak{m}=(\frac{\alpha-1}{2\alpha-1-\gamma})p$ and $\mathfrak{r}=\frac{d+\alpha}{2\alpha-1-\gamma}$,\\[2mm]
if $g:[0,+\infty[\times\mathbb{R}^d\longrightarrow \mathbb{R}$ is an external force such that we have 
$g \in \dot{\mathcal{W}}^{-\delta,\mathfrak{n},q}_\alpha([0,+\infty[\times\mathbb{R}^d)$, where $0<\delta<\alpha$ with $\mathfrak{n}=(\frac{\alpha-1}{3\alpha-1-\delta})p$ and $\mathfrak{s}= \frac{d+\alpha}{3\alpha-1-\delta}$.
\end{itemize}
If the quantity 
$$\|\vu_0\|_{T^{-\frac{\alpha}{p}}_{M^{p,q}}} + \|\theta_0\|_{T^{-\frac{\alpha}{\mathfrak{p}}}_{M^{\mathfrak{p},\mathfrak{q}}}} + \|\vf\|_{\dot{\mathcal{W}}^{-\gamma,\mathfrak{m}, \mathfrak{r}}_\alpha} + \|g\|_{\dot{\mathcal{W}}^{-\delta,\mathfrak{n},\mathfrak{s}}_\alpha},$$ 
is small enough, then there exists a \emph{global mild} solution $(\vu,\theta)$ of the fractional integral problem (\ref{Formulation_Integrale}) such that 
$$\vu \in \mathcal{M}_\alpha^{p,\frac{d+\alpha}{\alpha-1}}([0, +\infty[\times \mathbb{R}^d)\quad  \mbox{and} \quad \theta \in \mathcal{M}_\alpha^{\mathfrak{p},\frac{d+\alpha}{2\alpha-1}}([0, +\infty[\times \mathbb{R}^d).\\[3mm]$$
\end{Theoreme}
This result deserves several remarks. Indeed, we first note that, from the point of view of the homogeneity with respect to the dilations (\ref{Def_Scaling}), (\ref{Def_Scaling_Forces}) and (\ref{Def_Scaling_DonneeInitiales}), the spaces considered here are \emph{critical} since we have the identities 
$$\|\vu_{\lambda, 0}\|_{T^{-\frac{\alpha}{p}}_{M^{p,q}}}=\|\vu_{0}\|_{T^{-\frac{\alpha}{p}}_{M^{p,q}}} \quad \mbox{and}\quad \|\theta_{\lambda, 0}\|_{T^{-\frac{\alpha}{\mathfrak{p}}}_{M^{\mathfrak{p}, \mathfrak{q}}}}=\|\theta_{0}\|_{T^{-\frac{\alpha}{\mathfrak{p}}}_{M^{\mathfrak{p}, \mathfrak{q}}}}, \mbox{ for the initial data},$$
$$\|\vf_\lambda\|_{\dot{\mathcal{W}}^{-\gamma,\mathfrak{m}, \mathfrak{r}}_\alpha} = \|\vf\|_{\dot{\mathcal{W}}^{-\gamma,\mathfrak{m}, \mathfrak{r}}_\alpha} \quad \mbox{and}\quad \|g_\lambda\|_{\dot{\mathcal{W}}^{-\delta,\mathfrak{n},\mathfrak{s}}_\alpha}=\|g\|_{\dot{\mathcal{W}}^{-\delta,\mathfrak{n},\mathfrak{s}}_\alpha}, \mbox{ for the external forces},$$
$$\|\vu_\lambda\|_{\mathcal{M}_\alpha^{p,\frac{d+\alpha}{\alpha-1}}}=\|\vu\|_{\mathcal{M}_\alpha^{p,\frac{d+\alpha}{\alpha-1}}} \quad \mbox{and}\quad  \|\theta_\lambda\|_{\mathcal{M}_\alpha^{\mathfrak{p},\frac{d+\alpha}{2\alpha-1}}}=\|\theta\|_{\mathcal{M}_\alpha^{\mathfrak{p},\frac{d+\alpha}{2\alpha-1}}}, \mbox{ for the resolution spaces},$$
so we are working here in a \emph{complete} critical setting, and this is possible through the choice of the different indexes that define the functional spaces considered above. An interesting consequence of this fact is that, under a smallness condition over the initial data and the external forces, we can obtain mild solutions that are \emph{global} in time.\\ 

To the best of our knowledge this result is new for the Boussinesq fractional equations (\ref{EquationNSB_intro}) as it provides a quite general framework for mild solutions and we do not know if it is possible to consider a largest critical functional setting than the one presented here.\\

Let us make now some comments on the values of the parameters used here. Besides the dimension $d\geq 2$ and the fractional power of the Laplacian $1<\alpha<2$, the main parameter is the index $p$ that satisfies $4<\frac{3\alpha-2}{\alpha-1}<p\leq \frac{d+\alpha}{\alpha-1}$ and which determines the rest of the indexes of the functional spaces used for the initial data, the external forces and the resolution spaces.\\

In particular we remark that the upper bound $p\leq \frac{d+\alpha}{\alpha-1}$ is essential as it guarantees that the different Morrey spaces are well defined: indeed, this upper bound gives the conditions $p\leq q$, $\mathfrak{p}\leq\mathfrak{q}$, $\mathfrak{m}\leq\mathfrak{r}$ and $\mathfrak{n}\leq\mathfrak{s}$ for the Morrey spaces $M^{p,q}$ used for the initial data $\vu_0$, for the space $M^{\mathfrak{p}, \mathfrak{q}}$ used for the initial temperature $\theta_0$, for the spaces $\mathcal{M}_\alpha^{\mathfrak{m},\mathfrak{r}}$ and $\mathcal{M}_\alpha^{\mathfrak{n},\mathfrak{s}}$ used for the external forces $\vf$ and $g$, and for the spaces $\mathcal{M}_\alpha^{p, \frac{d+\alpha}{\alpha-1}}$ and $\mathcal{M}_\alpha^{\mathfrak{p}, \frac{d+\alpha}{2\alpha-1}}$ used as resolution spaces. \\

We note now that the lower bound $\frac{3\alpha-2}{\alpha-1}<p$ not only gives the conditions $1\leq p, \mathfrak{p}, \mathfrak{m}, \mathfrak{n}$ needed to correctly define the corresponding Morrey spaces, but it will also be crucial in order to perform some estimates over the bilinear quantity $B(\cdot, \cdot)$, specifically in the term $\displaystyle{\int_0^t} \mathfrak{p}_{t-s}\ast \mathrm{div}(\theta \vec{u})ds$ that takes into account the product $\theta \vu $: indeed, since the resolution spaces linked to each one of these variables are of different homogeneity, it is necessary to have enough integrability information to obtain the wished estimates. See Remark \ref{RemarqueValueP} below for more details on this particular point.\\

This lower bound for the parameter $p$ will also have other consequences. Indeed, parabolic Morrey spaces $\mathcal{M}^{p,q}_\alpha$ appear quite naturally in the study of mild solutions for the fractional Navier-Stokes equations as they provide a good functional framework for time-space estimates (see our recent work \cite{Chamorro_Mansais} as well as the computations performed below), moreover they also allow us to consider a very large class of initial data. \\

To make this point clear, let us for a moment focus our study in the functional spaces for the initial velocity field $\vu_0$. From a generic point of view, when dealing with homogeneous spaces -if in addition we ask a translation invariant property- it is known (see \cite{Meyer}) that Besov spaces $\dot{B}^{-\beta, \infty}_{\infty}(\mathbb{R}^d)$ are maximal in the sense that they contain all the functional spaces with the same homogeneity invariances. Recall that the homogeneous Besov spaces $\dot{B}^{-\beta, \infty}_\infty(\mathbb{R}^d)$ can be characterized by the condition
\begin{equation}\label{DefinitionBesov}
\|f\|_{\dot{B}^{-\beta, \infty}_\infty}=\underset{t>0}{\sup}\, t^{\frac{\beta}{\alpha}}\|\mathfrak{p}_t\ast f\|_{L^\infty}<+\infty,
\end{equation}
for any $f\in \mathcal{S}'/\mathcal{P}$. See \cite{Miao} and \cite{quad_non} for more details about this thermic characterization of Besov spaces. \\

However, if we consider the classical 3D Navier-Stokes equations (\emph{i.e.} when $\alpha=2$ and $\theta=0$ in (\ref{EquationNSB_intro})), it was shown in \cite{besov_ill_posed} that the largest homogeneous Besov space $\dot{B}^{-1, \infty}_{\infty}(\mathbb{R}^3)$, although critical with respect of the dilations associated to the structure of the Navier-Stokes equations, is not well suited to perform a fixed-point argument. As an alternative the space $BMO^{-1}(\mathbb{R}^3)\subset \dot{B}^{-1, \infty}_{\infty}(\mathbb{R}^3)$ was proposed in \cite{Koch}, which seems to be -at the best of our knowledge- the largest critical functional space for initial velocity fields $\vu_0$.\\ 

This situation changes completely if we consider the \emph{fractional} Navier-Stokes equations with $1<\alpha<2$, indeed, as it was pointed out in \cite{YuZhai} (see also \cite{Zhai} or \cite{Chamorro_Mansais}), an initial data\footnote{This Besov space corresponds with the scaling (\ref{Def_Scaling_DonneeInitiales}) for $\vu_0$.} $\vu_0\in \dot{B}^{-(\alpha-1), \infty}_\infty(\mathbb{R}^3)$ can be used without any problem in order to perform a fixed-point argument, in particular when we use as a resolution framework the critical parabolic Morrey space $\mathcal{M}^{p,\frac{d+\alpha}{\alpha-1}}_\alpha$ and this is due to the fact that we have the estimate 
\begin{equation}\label{ControlMorreyBesovData}
\|\mathfrak{p}_t\ast \vu_0\|_{\mathcal{M}^{p,\frac{d+\alpha}{\alpha-1}}_\alpha}\leq \|\vu_0\|_{\dot{B}^{-(\alpha-1), \infty}_\infty},
\end{equation}
as long as the parameter $p$ is \emph{small enough}. Indeed we have the following generic result which shows a useful relationship between the quantity $\mathfrak{p}_t\ast f$, considered in a parabolic Morrey space and the Besov norm of the function $f$:
\begin{Proposition}\label{Propo_Inclusion_MorreyParabolique_Besov}
Consider a function $f:\mathbb{R}^d\longrightarrow \mathbb{R}$ such that $f\in \dot{B}^{-\beta, \infty}_\infty(\mathbb{R}^d)$ for some parameter $0<\beta<d$. If $p,q$ are indexes such that 
 $1\leq p<\frac{\alpha}{\beta}$, $q=\frac{d+\alpha}{\beta}$ and $p\leq q$, then we have the following estimate:
$$\|\mathfrak{p}_t\ast f\|_{\mathcal{M}^{p,q}_\alpha}\leq C\|f\|_{\dot{B}^{-\beta, \infty}_\infty}.$$
\end{Proposition}
See a proof of this result in the appendix \ref{AppendixA}. In the general context of this article, this result states that, when dealing with a parabolic Morrey spaces $\mathcal{M}^{p,q}_\alpha$ as resolution spaces, we can work with an initial data $\vu_0$ (or $\theta_0$) that belongs to a suitable Besov space as long as the first parameter $p$ of the parabolic Morrey space is \emph{small}, which is unfortunately not the case of the Theorem \ref{Theorem_Triebel_Morrey} as we have the lower bound $\frac{3\alpha-2}{\alpha-1}<p$: the quantity $\|\mathfrak{p}_t\ast \vu_0\|_{\mathcal{M}_\alpha^{p, \frac{d+\alpha}{\alpha-1}}}$ can not be easily controlled by a Besov norm of the initial data $\vu_0$ (in the sense of the estimate (\ref{ControlMorreyBesovData}) above) as the arguments displayed in the proof of the Proposition \ref{Propo_Inclusion_MorreyParabolique_Besov} fail in this case (the constraint $\frac{3\alpha-2}{\alpha-1}<p$ is incompatible with the condition $p<\frac{\alpha}{\beta}=\frac{\alpha}{\alpha-1}$ if $1<\alpha<2$). This particular point forced us to consider for the initial data $\vu_0$ the Triebel-Lizorkin-Morrey spaces introduced in the expression (\ref{Def_EspaceTriebelLizorkinMorrey}). Note that the context is even more delicate if we want to consider an initial data $\theta_0\in \dot{B}^{-(2\alpha-1),\infty}_\infty(\mathbb{R}^d)$ since in this case we have $\beta=2\alpha-1>\alpha$ which would give the condition $p<\frac{\alpha}{\beta}=\frac{\alpha}{2\alpha-1}<1$ and the previous result can not be applied. \\

In the following result, we show that, as long as we want to work with the critical parabolic Morrey spaces $\mathcal{M}_\alpha^{p,\frac{d+\alpha}{\alpha-1}}$ as resolution spaces, then the Triebel-Lizorkin-Morrey spaces introduced in the expression (\ref{Def_EspaceTriebelLizorkinMorrey}) above are quite natural:
\begin{Theoreme}\label{Theorem_Equivalences}
For $1<\alpha<2$, if $(\vu,\theta)$ is a mild solution of the fractional integral problem (\ref{Formulation_Integrale}) such that 
$$\vu \in \mathcal{M}_\alpha^{p,\frac{d+\alpha}{\alpha-1}}([0, +\infty[\times \mathbb{R}^d)\quad  \mbox{and} \quad \theta \in \mathcal{M}_\alpha^{\mathfrak{p},\frac{d+\alpha}{2\alpha-1}}([0, +\infty[\times \mathbb{R}^d),$$
where the parameter $p$ satisfies the conditions $\frac{3\alpha-2}{\alpha-1}<p\leq \frac{d+\alpha}{\alpha-1}$, then the initial velocity field $\vu_0$ belongs to the space $T^{-\frac{\alpha}{p}}_{M^{p,q}}(\mathbb{R}^d)$ with $\frac{1}{q}=\frac{\alpha-1}{d+\alpha}-\frac{\alpha}{d}(\frac{1}{p}-\frac{\alpha-1}{d+\alpha})$ and the initial temperature $\theta_0$ belongs to the space  $T^{-\frac{\alpha}{\mathfrak{p}}}_{M^{\mathfrak{p},\mathfrak{q}}}(\mathbb{R}^d)$ with $\mathfrak{p}=(\frac{\alpha-1}{2\alpha-1})p$ and $\frac{1}{\mathfrak{q}}=(\frac{2\alpha-1}{d+\alpha})-\frac{\alpha}{d}(\frac{1}{\mathfrak{p}}-\frac{2\alpha-1}{d+\alpha})$.
\end{Theoreme}
In order to establish this theorem we will prove the following equivalences
$$\|\mathfrak{p}_t\ast \vu_0\|_{ \mathcal{M}_\alpha^{p,\frac{d+\alpha}{\alpha-1}}}\simeq \|\vu_0\|_{T^{-\frac{\alpha}{p}}_{M^{p,q}}} \quad \mbox{and}\quad \|\mathfrak{p}_t\ast \theta_0\|_{\mathcal{M}_\alpha^{\mathfrak{p},\frac{d+\alpha}{2\alpha-1}}}\simeq \|\theta_0\|_{T^{-\frac{\alpha}{\mathfrak{p}}}_{M^{\mathfrak{p},\mathfrak{q}}}},$$
which show how the Triebel-Lizorkin-Morrey spaces appear naturally in this context.\\

To the best of our knowledge, Theorems \ref{Theorem_Triebel_Morrey} and \ref{Theorem_Equivalences} are new in the context of the fractional Boussinesq equations and although critical parabolic Morrey spaces can be used to perform a fixed point argument, it seems that due to the structural homogeneity structure of the system, the homogeneous Besov spaces can not be considered for the initial data in this framework: the possibility of considering mild solutions in a critical setting with an initial data $\vu_0\in \dot{B}^{-(\alpha-1),\infty}_\infty(\mathbb{R}^d)$ and $\theta_0\in \dot{B}^{-(2\alpha-1),\infty}_\infty(\mathbb{R}^d)$ seems to be a new open problem.\\

The plan of the article is the following. In Section \ref{Secc_Notations} we recall some notation and useful estimates and properties of the functional spaces involved in this article, in Section \ref{Secc_Proof_Theorem_Triebel_Morrey} we give the proof of the Theorem \ref{Theorem_Triebel_Morrey} and in Section \ref{Secc_Proof_Theorem_Equivalences} we give the proof the Theorem \ref{Theorem_Equivalences}. 
\section{Kernel estimates, parabolic Riesz potentials and some inequalities}\label{Secc_Notations}
In the section we gather some essential facts that will be used in the proofs of our theorems.
\begin{itemize}
\item {\bf Pointwise estimates for the fractional heat kernel.} One of the essential features that will be displayed on this article is related to some specific pointwise estimates for the fractional heat kernel $\mathfrak{p}_t$. Indeed, we first introduce the following notation: let $\sigma_\rho \in \mathcal{C}^{\infty}(\mathbb{R}^d \setminus \{0\})$ be a positive, homogeneous function on $\mathbb{R}^d\setminus \{0\}$, with homogeneous degree $\rho>-d$ (\emph{i.e.} we have $\sigma_\rho(\lambda \cdot)=\lambda^{\rho}\sigma_\rho(\cdot)$ for all $\lambda>0$). Then, for $1<\alpha<2$ and for $t>0$, we define the function $\mathcal{K}^\rho_t:\mathbb{R}^d\longrightarrow \mathbb{R}$ in the Fourier level by the condition 
$$\widehat{\mathcal{K}^\rho_t}(\xi)=\sigma_\rho(\xi)e^{-t|\xi|^\alpha}.$$
It is worth mentioning that the condition $\rho>-d$ is only used to ensure that $\sigma_\rho(\xi)e^{-t|\xi|^\alpha}$ is integrable with respect to $\xi$ in a neighborhood of the origin. With this notation at hand we can present now a pointwise estimate for such kernels, which will immediately grant us useful integrability properties: indeed, for $1<\alpha<2$, for $t>0$ and for some constant $C>0$, we have:
\begin{equation}\label{Prop_estim_ptws}
|\mathcal{K}^{\rho}_t(x)| \leq \frac{C}{(t^\frac{1}{\alpha}+|x|)^{d+\rho}}.
\end{equation}
A proof of this estimate can be consulted in \cite{quad_non}. From this pointwise inequality we easily deduce the $L^p$ control if $\rho>-\frac{d}{p'}$
\begin{equation}\label{Prop_estim_Lp}
\|\mathcal{K}^{\rho}_t\|_{L^p}\leq C t^{-\frac{\rho+d(1-\frac{1}{p})}{\alpha}},
\end{equation}
valid for $1\leq p\leq +\infty$. For more details and properties of this type of semi-groups, see the book \cite{kolokoltsov}.
\item{\bf Parabolic Riesz potential.} For a locally integrable function $f:\mathbb{R}\times \mathbb{R}^d\longrightarrow \mathbb{R}$, we define the parabolic Riesz potential $\mathrm{I}_\beta$ with $0<\beta<d+\alpha$ by the expression
\begin{equation}\label{Def_ParabolicRiesz}
\mathrm{I}_\beta(f)(t,x) = \displaystyle{\int_{\mathbb{R}\times \mathbb{R}^d}} \frac{1}{(|t-s|^\frac{1}{\alpha}+|x-y|)^{d+\alpha-\beta}}f(s,y)dyds.
\end{equation}
As we shall see, the parabolic Riesz potential $\mathrm{I}_\beta$ will naturally appear as a corollary of the estimates (\ref{Prop_estim_ptws}) given previously (for this it will be ``enough'' to set $\rho=\alpha-\beta$ in the kernel inequality above). 
\item{\bf Parabolic Morrey spaces inequalities.} We recall here two important inequalities in the context of the parabolic Morrey spaces presented in (\ref{Definition_MorreyPara}). The first result is nothing but a suitable version of the H\"older inequalities: 
\begin{Proposition}
Let $1<\alpha<2$ be fixed. Let also $p_1,p_2,q_1,q_2$ denote real parameters such that     $1\leq p_1\leq q_1 <+\infty$ and $1\leq p_2\leq q_2 <+\infty$. If $f\in \mathcal{M}_\alpha^{p_1,q_1}(\mathbb{R}\times \mathbb{R}^d)$ and if $g\in \mathcal{M}_\alpha^{p_2,q_2}(\mathbb{R}\times \mathbb{R}^d)$, then we have the following version of the H\"older inequality in the parabolic Morrey space setting:
$$\|fg\|_{\mathcal{M}_\alpha^{p,q}} \leq \|f\|_{\mathcal{M}_\alpha^{p_1,q_1}} \|g\|_{\mathcal{M}_\alpha^{p_2,q_2}},$$
where $\frac{1}{p}=\frac{1}{p_1}+\frac{1}{p_2}$ and $\frac{1}{q}=\frac{1}{q_1}+\frac{1}{q_2}$.
\end{Proposition}
The second result present some boundedness of the parabolic Riesz potentials in parabolic Morrey spaces: 
\begin{Proposition}\label{Prop_continuity_Riesz}
Let $1<p\leq q <+\infty$ and $1<\alpha<2$. If $f$ is a measurable function on $\mathbb{R}\times \mathbb{R}^d$ such that $f\in \mathcal{M}_\alpha^{p,q}(\mathbb{R}\times \mathbb{R}^d)$ and if  $0<\beta<\frac{d+\alpha}{q}$, then we have the following estimate
$$\left\|\mathrm{I}_\beta(f)\right\|_{\mathcal{M}_\alpha^{\frac{p}{\lambda},\frac{q}{\lambda}}} \leq C \|f\|_{\mathcal{M}_\alpha^{p,q}},$$
where $\lambda=1-\beta\frac{q}{d+\alpha}$.
\end{Proposition}
See a proof of this fact in \cite{quad_non}.
\end{itemize}
\section{Proof of the Theorem \ref{Theorem_Triebel_Morrey}}\label{Secc_Proof_Theorem_Triebel_Morrey}
We first remark that the integral formulation (\ref{Formulation_Integrale}) can be rewritten in the following manner
\begin{equation}\label{Equation_PointFixe}
\vU=\vU_0+\vec{F}+L(\vU)+B(\vU,\vU),
\end{equation}
where $\vU=\left(\begin{array}{c}\vu\\ \theta\end{array}\right)$ is a $3+1$ vector, $\vU_0 = \left(\begin{array}{c}
\mathfrak{p}_t\ast \vec{u_0} \\\mathfrak{p}_t\ast \theta_0\end{array}\right)$ is the initial data, $\vec{F} = \left(\begin{array}{c} \displaystyle{\int_0^t} \mathfrak{p}_{t-s}\ast\mathbb{P}(\vec{f})ds\\
 \displaystyle{\int_0^t} \mathfrak{p}_{t-s}\ast g\,ds\end{array}\right)$ is the external force (recall that $\vf$ and $g$ are given), $L(\vU)=\left(\begin{array}{c} \displaystyle{\int_0^t} \mathfrak{p}_{t-s}\ast\mathbb{P}(\theta \vec{e}_d)ds\\  0\end{array}\right)$ is a linear term and $B(\vU,\vU)=\left(\begin{array}{c} -\displaystyle{\int_0^t} \mathfrak{p}_{t-s}\ast\mathbb{P}\left(\mathrm{div}(\vec{u}\otimes \vec{u})\right)ds\\     \displaystyle{\int_0^t} \mathfrak{p}_{t-s}\ast \mathrm{div}(\theta \vec{u})ds\end{array}\right)$ is a bilinear term.\\ 
 
\noindent The equations of the form (\ref{Equation_PointFixe}) can be easily studied as we have the following version of the Banach fixed point theorem: 
\begin{Lemme}\label{Lemme_Banach_Picard}
Let $(E,\|\cdot\|_E)$ be a Banach space. Consider $\vU_0 \in E$ and let $L:E \longrightarrow E$ be a linear operator and $B:E\times E \longrightarrow E$ be a bilinear operator that satisfy the estimates
\begin{eqnarray*}
\|L(\vU)\|_E &\leq &C_L \|\vU\|_E,\\[2mm]
\|B(\vU,\vec{V})\|_E &\leq &C_B \|\vU\|_E \|\vec{V}\|_E, 
\end{eqnarray*}
 for all $\vU, \vec{V}\in E$. Assume moreover that $C_L<\frac{1}{3}$. If  we have $\|\vU_0\|_E \leq \delta/2$ and $\|\vec{F}\|_E \leq \delta/2$ for a constant $\delta>0$ such that $\delta < \frac{1}{9C_B}$, then there exists a unique $\vU \in E$ such that $\|\vU\|_E\leq 3 \delta$ which is a solution of the equation (\ref{Equation_PointFixe}).
\end{Lemme}
See \cite[Théorème 7.1.1]{Chamorro_Livre} for a proof of this lemma.\\

In what follows we will apply this general result to construct mild solutions of the integral problem (\ref{Formulation_Integrale}) and of course the choice of the resolution space $(E,\|\cdot\|_E)$ will condition the nature of the results that we want to obtain. Indeed, we are interested here in considering a resolution space that is \emph{critical} with respect to the scaling (\ref{Def_Scaling}), however, since the scaling for the variable $\vu$ is different from the one of the variable $\theta$, we will separate the information in the following manner: for the vector $\vU=\left(\begin{array}{c}\vu\\ \theta\end{array}\right)$ we will consider the norm 
\begin{equation}\label{Def_NormeGlobale}
\|\vU\|_{E}=\left\|\left(\begin{array}{c}\vu\\ \theta\end{array}\right)\right\|_{E}=\|\vu\|_{E_1}+\mathfrak{C}\|\theta\|_{E_2},
\end{equation}
where $(E_1, \|\cdot\|_{E_1})$ and $(E_2, \|\cdot\|_{E_2})$ are two Banach function spaces and where $\mathfrak{C}=\mathfrak{C}(\alpha,d)\gg 1$ is a constant that will be fixed later. Thus, if we define for some $\lambda>0$ the vector $\vU_\lambda=\left(\begin{array}{c}\vu_\lambda\\ \theta_\lambda\end{array}\right)$, where the scaled functions $\vu_\lambda$ and $\theta_\lambda$ are given in the expression (\ref{Def_Scaling}) above, we will say that the norm $\|\cdot\|_E$ is \emph{critical} with respect to the scaling (\ref{Def_Scaling}) if we have the identity 
\begin{eqnarray*}
\|\vU_\lambda\|_{E}&=&\|\vu_\lambda\|_{E_1}+\mathfrak{C}\|\theta_\lambda\|_{E_2}\\
&=&\|\vu\|_{E_1}+\mathfrak{C}\|\theta\|_{E_2}=\|\vU\|_{E}.
\end{eqnarray*}
Note that the homogeneity property asked for the Banach spaces $E_1$ and $E_2$ is different since the scaling of the variables $\vu$ and $\theta$ is different.\\ 

As we will see, one very important feature of \emph{critical} norms is that it is possible to obtain mild solutions that are \emph{global} in time (as they preserve the scaling structure of the equation). However, in order to work in a fully \emph{critical} setting, we also need to choose a suitable framework for the initial data $(\vu_0, \theta_0)$ and for the external forces $(\vf, g)$. We will thus say that a Banach space $(\mathcal{E}, \|\cdot\|_\mathcal{E})$ is a critical space for the initial data $\left(\begin{array}{c}\vu_0\\ \theta_0\end{array}\right)$ if we have the identity
\begin{equation}\label{Def_NormeGlobale_DonneesInitiales}
\left\|\left(\begin{array}{c}\vu_{\lambda,0}\\ \theta_{\lambda,0}\end{array}\right)\right\|_{\mathcal{E}}=\|\vu_{\lambda,0}\|_{\mathcal{E}_1}+\|\theta_{\lambda,0}\|_{\mathcal{E}_2}=\|\vu_{0}\|_{\mathcal{E}_1}+\|\theta_{0}\|_{\mathcal{E}_2}=\left\|\left(\begin{array}{c}\vu_{0}\\ \theta_{0}\end{array}\right)\right\|_{\mathcal{E}},
\end{equation}
where the scaled functions $\vu_{\lambda,0}$ and $\theta_{\lambda, 0}$ are given in the expression (\ref{Def_Scaling_DonneeInitiales}) above. In the same spirit, we will say that a Banach space $(\mathcal{F}, \|\cdot\|_\mathcal{F})$ is a critical space for the external forces  $\left(\begin{array}{c}\vf\\ g\end{array}\right)$ if we have the identity 
\begin{equation}\label{Def_NormeGlobale_ForcesExterieurs}
\left\|\left(\begin{array}{c}\vf_{\lambda}\\ g_{\lambda}\end{array}\right)\right\|_{\mathcal{F}}=\|\vf_{\lambda}\|_{\mathcal{F}_1}+\|g_{\lambda}\|_{\mathcal{F}_2}=\|\vf\|_{\mathcal{F}_1}+\|g\|_{\mathcal{F}_2}=\left\|\left(\begin{array}{c}\vf\\ g\end{array}\right)\right\|_{\mathcal{F}},
\end{equation}
where the scaled functions $\vf_\lambda$ and $g_\lambda$ are defined in the formula (\ref{Def_Scaling_Forces}) above.\\

With this generic critical framework at hand, we will construct in this article global mild solutions for the system (\ref{Formulation_Integrale}) and in order to apply the Lemma \ref{Lemme_Banach_Picard}, once we have fixed a critical function Banach space $(E, \|\cdot\|_E)$ as a resolution space, we will need to establish the following estimates (where we are considering the formulation (\ref{Equation_PointFixe}) above):
\begin{itemize}
\item For the initial data $\vU_0$: 
\begin{eqnarray}
\|\vU_0\|_{E} = \left\|\left(\begin{array}{c}
\mathfrak{p}_t\ast \vec{u_0} \\\mathfrak{p}_t\ast \theta_0\end{array}\right)\right\|_E&= &\|\mathfrak{p}_t\ast \vec{u_0}\|_{E_1}+\mathfrak{C}\|\mathfrak{p}_t\ast \theta_0\|_{E_2}\notag\\
&\leq &C\|\vu_0\|_{\mathcal{E}_1}+C\mathfrak{C}\|\theta_0\|_{\mathcal{E}_2}\label{Estimations_Generales_DonneesInitiales}
\end{eqnarray}
\item For the external force $\vec{F}$: 
\begin{eqnarray}
\|\vec{F}\|_{E} = \left\|\left(\begin{array}{c} \displaystyle{\int_0^t} \mathfrak{p}_{t-s}\ast\mathbb{P}(\vec{f})ds\\
 \displaystyle{\int_0^t} \mathfrak{p}_{t-s}\ast g\,ds\end{array}\right)\right\|_E&= &\left\| \displaystyle{\int_0^t} \mathfrak{p}_{t-s}\ast\mathbb{P}(\vec{f})ds\right\|_{E_1}+\mathfrak{C}\left\| \displaystyle{\int_0^t} \mathfrak{p}_{t-s}\ast g\,ds\right\|_{E_2}\notag\\
&\leq &C\|\vf\|_{\mathcal{F}_1}+C\mathfrak{C}\|g\|_{\mathcal{F}_2}\label{Estimations_Generales_ForcesExterieures}
\end{eqnarray}

\item For the linear term $L(\vU)$:
\begin{eqnarray*}
\|L(\vU)\|_{E}= \left\|\left(\begin{array}{c} \displaystyle{\int_0^t} \mathfrak{p}_{t-s}\ast\mathbb{P}(\theta \vec{e}_d)ds\\  0\end{array}\right)
\right\|_E&= &\left\| \displaystyle{\int_0^t} \mathfrak{p}_{t-s}\ast\mathbb{P}(\theta \vec{e}_d)ds\right\|_{E_1}\leq C_0\|\theta\|_{E_2},
\end{eqnarray*}
recalling that we have $\|\theta\|_{E_2}\leq \frac{1}{\mathfrak{C}}(\|\vu\|_{E_1}+\mathfrak{C}\|\theta\|_{E_2})$, we can write
\begin{eqnarray}
\|L(\vU)\|_{E}&\leq &\frac{C_0}{\mathfrak{C}}(\|\vu\|_{E_1}+\mathfrak{C}\|\theta\|_{E_2})=\frac{C_0}{\mathfrak{C}}\|\vU\|_E\notag\\
&\leq & C_L\|\vU\|_E,\label{Estimations_Generales_TermeLineaire}
\end{eqnarray}
where we have $C_L=\frac{C_0}{\mathfrak{C}}$.
\item For the bilinear term $B(\vU, \vU)$:
\begin{eqnarray*}
\|B(\vU, \vU)\|_{E}&=&\left\|\left(\begin{array}{c} -\displaystyle{\int_0^t} \mathfrak{p}_{t-s}\ast\mathbb{P}\left(\mathrm{div}(\vec{u}\otimes \vec{u})\right)ds\\   \displaystyle{\int_0^t} \mathfrak{p}_{t-s}\ast \mathrm{div}(\theta \vec{u})ds\end{array}\right)\right\|_E\\
&\leq &\left\| \displaystyle{\int_0^t} \mathfrak{p}_{t-s}\ast\mathbb{P}\left(\mathrm{div}(\vec{u}\otimes \vec{u})\right)ds\right\|_{E_1}+\mathfrak{C}\left\|  \displaystyle{\int_0^t} \mathfrak{p}_{t-s}\ast \mathrm{div}(\theta \vec{u})ds\right\|_{E_2}\\
&\leq & C_1\|\vu\|_{E_1}\|\vu\|_{E_1}+C_2\mathfrak{C}\|\vu\|_{E_1}\|\theta\|_{E_2},
\end{eqnarray*}
and since we have 
$$\|\vu\|_{E_1}\leq\|\vu\|_{E_1}+\mathfrak{C}\|\theta\|_{E_2}=\|\vU\|_{E}\quad \mbox{and}\quad \mathfrak{C}\|\theta\|_{E_2}\leq \|\vu\|_{E_1}+\mathfrak{C}\|\theta\|_{E_2}=\|\vU\|_{E},$$
we can write 
\begin{equation}
\|B(\vU, \vU)\|_{E}\leq (C_1+C_2)\|\vU\|_{E}\|\vU\|_{E}.\label{Estimations_Generales_TermeBilineaire}
\end{equation}
\end{itemize}
Several remarks are in order here. First note that the only purpose of the constant $\mathfrak{C}$ introduced in the norm (\ref{Def_NormeGlobale}) is to make the constant $C_L$ in the estimate (\ref{Estimations_Generales_TermeLineaire}) small enough, which is necessary to apply the Lemma \ref{Lemme_Banach_Picard}. Remark next that since we want to work in a critical setting, the constants that appear in all the previous estimates will not depend in the time variable: thus, under a ``smallness'' hypothesis (stated in terms of suitable norms) over the initial data $(\vu_0, \theta_0)$ and over the external forces $(\vf, g)$, we will obtain \emph{global in time} mild solutions for the problem (\ref{Formulation_Integrale}). This \emph{smallness} condition is then expressed by the estimate 
$$C\|\vu_0\|_{\mathcal{E}_1}+C\mathfrak{C}\|\theta_0\|_{\mathcal{E}_2}+C\|\vf\|_{\mathcal{F}_1}+C\mathfrak{C}\|g\|_{\mathcal{F}_2}<\frac{1}{9(C_1+ C_2)}.$$
Note now that, following the structure of the norm $\|\cdot\|_E$ of the resolution space given in (\ref{Def_NormeGlobale}), the norm naturally ``linked'' to the variable $\vu$ is $\|\cdot\|_{E_1}$ while the norm ``linked'' to the variable $\theta$ is $\|\cdot\|_{E_2}$ (this peculiarity is also reflected in the initial data as well as in the external forces, see the expressions (\ref{Def_NormeGlobale_DonneesInitiales}) and (\ref{Def_NormeGlobale_ForcesExterieurs}) above) and most of the estimates presented in (\ref{Estimations_Generales_DonneesInitiales})-(\ref{Estimations_Generales_TermeBilineaire}) respect this separation of norms. However,  since the Boussinesq system contains coupled terms, the study of the linear term $L(\vU)$ as well as the bilinear term $B(\vU, \vU)$ will \emph{structurally} require some mixed norm estimates (see the estimates (\ref{Estimations_Generales_TermeLineaire}) and (\ref{Estimations_Generales_TermeBilineaire}) above), and this will demand a certain flexibility in the norms considered.\\

We will now establish the general estimates (\ref{Estimations_Generales_DonneesInitiales})-(\ref{Estimations_Generales_TermeBilineaire}) within the functional framework stated in Theorem \ref{Theorem_Triebel_Morrey}.

\begin{itemize}
\item For the initial data $\vU_0=\left(\begin{array}{c}
\mathfrak{p}_t\ast \vec{u_0} \\\mathfrak{p}_t\ast \theta_0\end{array}\right)$ we write

$$\|\vU_0\|_{E} = \left\|\left(\begin{array}{c}
\mathfrak{p}_t\ast \vec{u_0} \\\mathfrak{p}_t\ast \theta_0\end{array}\right)\right\|_E=\underbrace{\|\mathfrak{p}_t\ast \vec{u_0}\|_{\mathcal{M}_\alpha^{p,\frac{d+\alpha}{\alpha-1}}}}_{(1_D)}+\underbrace{\mathfrak{C}\|\mathfrak{p}_t\ast \theta_0\|_{\mathcal{M}_\alpha^{\mathfrak{p},\frac{d+\alpha}{2\alpha-1}}}}_{(2_D)}.$$
For the first norm $(1_D)$ in the right-hand side, using the definition of the parabolic Morrey norm given in (\ref{Definition_MorreyPara}), we have
$$\|\mathfrak{p}_t\ast \vu_0\|_{\mathcal{M}_\alpha^{p,\frac{d+\alpha}{\alpha-1}}}=\underset{r>0}{\mathrm{sup}} \ \underset{(t,x)\in [0,+\infty[ \times \mathbb{R}^d}{\mathrm{sup}} \ \frac{1}{r^{(d+\alpha)(\frac{1}{p}-\frac{\alpha-1}{d+\alpha})}} \left(\displaystyle{\iint_{\{|t-s|^\frac{1}{\alpha}+|x-y|<r\}}}|\mathfrak{p}_s\ast \vu_0(y)|^pdyds \right)^\frac{1}{p},$$
 Since the integration domain $\{|t-s|^\frac{1}{\alpha}+|x-y|<r\}$ is included in the set $[0,+\infty[\times B(x,r)$ we find the upper bound
$$\|\mathfrak{p}_t\ast \vu_0\|_{\mathcal{M}_\alpha^{p,\frac{d+\alpha}{\alpha-1}}}\leq \underset{r>0}{\mathrm{sup}} \ \underset{x\in \mathbb{R}^d}{\mathrm{sup}} \ \frac{1}{r^{(d+\alpha)(\frac{1}{p}-\frac{\alpha-1}{d+\alpha})}} \left(\displaystyle{\int_{B(x,r)}\int_{0}^{+\infty}}|\mathfrak{p}_s\ast \vu_0(y)|^pds dy\right)^\frac{1}{p}.$$
We thus have
\begin{eqnarray*}
\|\mathfrak{p}_t\ast \vu_0\|_{\mathcal{M}_\alpha^{p,\frac{d+\alpha}{\alpha-1}}}\leq \underset{r>0}{\mathrm{sup}} \ \underset{x\in \mathbb{R}^d}{\mathrm{sup}} \ \frac{1}{r^{(d+\alpha)(\frac{1}{p}-\frac{\alpha-1}{d+\alpha})}} \left(\int_{B(x,r)}\left[\left(\int_{0}^{+\infty}|\mathfrak{p}_s\ast \vu_0(y)|^p ds\right)^{\frac{1}{p}}\right]^p dy\right)^\frac{1}{p}\\
\leq  \underset{r>0}{\mathrm{sup}} \ \underset{x\in \mathbb{R}^d}{\mathrm{sup}} \ \frac{1}{r^{(d+\alpha)(\frac{1}{p}-\frac{\alpha-1}{d+\alpha})}} \frac{r^{d(\frac{1}{p}-\frac{1}{q})}}{r^{d(\frac{1}{p}-\frac{1}{q})}}\left\|\left(\int_{0}^{+\infty}|\mathfrak{p}_s\ast \vu_0(y)|^pds\right)^{\frac{1}{p}}\right\|_{L^p(B(x,r))}.
\end{eqnarray*}
Since by hypothesis we have $\frac{1}{q}=\frac{\alpha-1}{d+\alpha}-\frac{\alpha}{d}(\frac{1}{p}-\frac{\alpha-1}{d+\alpha})$, then $(d+\alpha)(\frac{1}{p}-\frac{\alpha-1}{d+\alpha})=d(\frac{1}{p}-\frac{1}{q})$ and we obtain
\begin{eqnarray*}
\|\mathfrak{p}_t\ast \vu_0\|_{\mathcal{M}_\alpha^{p,\frac{d+\alpha}{\alpha-1}}}\leq  \underset{r>0}{\mathrm{sup}} \ \underset{x\in \mathbb{R}^d}{\mathrm{sup}} \  \frac{1}{r^{d(\frac{1}{p}-\frac{1}{q})}}\left\|\left(\int_{0}^{+\infty}|\mathfrak{p}_s\ast \vu_0(y)|^pds\right)^{\frac{1}{p}}\right\|_{L^p(B(x,r))}=\|\vu_0\|_{T^{-\frac{\alpha}{p}}_{M^{p,q}}}.
\end{eqnarray*}
For the second norm $(2_D)$ we have in a similar manner
\begin{eqnarray*}
\|\mathfrak{p}_t\ast \theta_0\|_{\mathcal{M}_\alpha^{\mathfrak{p},\frac{d+\alpha}{2\alpha-1}}}&=&\underset{r>0}{\mathrm{sup}} \ \underset{(t,x)\in [0,+\infty[ \times \mathbb{R}^d}{\mathrm{sup}} \ \frac{1}{r^{(d+\alpha)(\frac{1}{\mathfrak{p}}-\frac{2\alpha-1}{d+\alpha})}} \left(\displaystyle{\iint_{\{|t-s|^\frac{1}{\alpha}+|x-y|<r\}}}|\mathfrak{p}_s\ast \theta_0(y)|^\mathfrak{p} dyds \right)^\frac{1}{\mathfrak{p}}\\
&\leq & \underset{r>0}{\mathrm{sup}} \ \underset{x\in \mathbb{R}^d}{\mathrm{sup}} \ \frac{1}{r^{(d+\alpha)(\frac{1}{\mathfrak{p}}-\frac{2\alpha-1}{d+\alpha})}} \left(\displaystyle{\int_{B(x,r)}\int_{0}^{+\infty}}|\mathfrak{p}_s\ast \theta_0(y)|^\mathfrak{p}ds dy\right)^\frac{1}{\mathfrak{p}},
\end{eqnarray*}
and the previous expression can be rewritten as
$$\|\mathfrak{p}_t\ast \theta_0\|_{\mathcal{M}_\alpha^{\mathfrak{p},\frac{d+\alpha}{2\alpha-1}}}\leq   \underset{r>0}{\mathrm{sup}} \ \underset{x\in \mathbb{R}^d}{\mathrm{sup}} \ \frac{1}{r^{(d+\alpha)(\frac{1}{\mathfrak{p}}-\frac{2\alpha-1}{d+\alpha})}} \frac{r^{d(\frac{1}{\mathfrak{p}}-\frac{1}{\mathfrak{q}})}}{r^{d(\frac{1}{\mathfrak{p}}-\frac{1}{\mathfrak{q}})}}\left\|\left(\int_{0}^{+\infty}|\mathfrak{p}_s\ast \theta_0(y)|^\mathfrak{p}ds\right)^{\frac{1}{\mathfrak{p}}}\right\|_{L^\mathfrak{p}(B(x,r))}.$$
As by hypothesis we have $(d+\alpha)(\frac{1}{\mathfrak{p}}-\frac{2\alpha-1}{d+\alpha})=d(\frac{1}{\mathfrak{p}}-\frac{1}{\mathfrak{q}})$ (indeed, recall that we have the identity $\frac{1}{\mathfrak{q}}=\frac{2\alpha-1}{d+\alpha}-\frac{\alpha}{d}(\frac{1}{\mathfrak{p}}-\frac{2\alpha-1}{d+\alpha})$), we thus obtain 
$$\|\mathfrak{p}_t\ast \theta_0\|_{\mathcal{M}_\alpha^{\mathfrak{p},\frac{d+\alpha}{2\alpha-1}}}\leq   \underset{r>0}{\mathrm{sup}} \ \underset{x\in \mathbb{R}^d}{\mathrm{sup}} \   \frac{1}{r^{d(\frac{1}{\mathfrak{p}}-\frac{1}{\mathfrak{q}})}}\left\|\left(\int_{0}^{+\infty}|\mathfrak{p}_s\ast \theta_0(y)|^\mathfrak{p}ds\right)^{\frac{1}{\mathfrak{p}}}\right\|_{L^\mathfrak{p}(B(x,r))}=\|\theta_0\|_{T^{\frac{-\alpha}{\mathfrak{p}}}_{M^{\mathfrak{p, \mathfrak{q}}}}}.$$
We can finally write, with the two previous estimates for $(1_D)$ and $(2_D)$:
$$\|\vU_0\|_{E} = \|\mathfrak{p}_t\ast \vec{u_0}\|_{\mathcal{M}_\alpha^{p,\frac{d+\alpha}{\alpha-1}}}+\mathfrak{C}\|\mathfrak{p}_t\ast \theta_0\|_{\mathcal{M}_\alpha^{\mathfrak{p},\frac{d+\alpha}{2\alpha-1}}}\leq \|\vu_0\|_{T^{-\frac{\alpha}{p}}_{M^{p,q}}}+ \mathfrak{C}\|\theta_0\|_{T^{\frac{-\alpha}{\mathfrak{p}}}_{M^{\mathfrak{p, \mathfrak{q}}}}}<+\infty.$$

\item We study the linear term $L(\vU)$ and we want to establish the estimate (\ref{Estimations_Generales_TermeLineaire}). For this we write
$$\|L(\vU)\|_{E} = \left\|\left(\begin{array}{c} \displaystyle{\int_0^t} \mathfrak{p}_{t-s}\ast\mathbb{P}(\theta \vec{e}_d)ds\\  0\end{array}\right)\right\|_E\leq \left\| \displaystyle{\int_0^t} \mathfrak{p}_{t-s}\ast\mathbb{P}(\theta \vec{e}_d)ds\right\|_{\mathcal{M}_\alpha^{p,\frac{d+\alpha}{\alpha-1}}}.$$
Now, since the Leray projector $\mathbb{P}$ is a vector of pseudo-differential operators with homogeneous symbol of degree $0$, smooth outside the origin, using the estimate (\ref{Prop_estim_ptws}) with $\rho=0$, we have the pointwise estimate (using the definition of the Riesz potentials (\ref{Def_ParabolicRiesz})):
$$\left| \displaystyle{\int_0^t} \mathfrak{p}_{t-s}\ast\mathbb{P}(\theta \vec{e}_d)ds\right| \leq C\displaystyle{\int_{\mathbb{R}\times \mathbb{R}^d}} \frac{1}{(|t-s|^\frac{1}{\alpha}+|x-y|)^{d}}|\theta(s,y)|dyds= C\,\mathrm{I}_{\alpha}(|\theta|)(t,x),$$
from which we easily derive the control 
$$\left\| \displaystyle{\int_0^t} \mathfrak{p}_{t-s}\ast\mathbb{P}(\theta \vec{e}_d)ds\right\|_{\mathcal{M}_\alpha^{p,\frac{d+\alpha}{\alpha-1}}}\leq C \|\mathrm{I}_{\alpha}(|\theta|) \|_{\mathcal{M}_\alpha^{p,\frac{d+\alpha}{\alpha-1}}}.$$
Applying Proposition \ref{Prop_continuity_Riesz} with $\lambda=\frac{\alpha-1}{2\alpha-1}$, we can write
$$\|\mathrm{I}_{\alpha}(|\theta|)\|_{\mathcal{M}_\alpha^{p,\frac{d+\alpha}{\alpha-1}}} \leq C \|\theta\|_{\mathcal{M}_\alpha^{(\frac{\alpha-1}{2\alpha-1})p,\frac{d+\alpha}{2\alpha-1}}}= C \|\theta\|_{\mathcal{M}_\alpha^{\mathfrak{p},\frac{d+\alpha}{2\alpha-1}}},$$
(recall that $\mathfrak{p}=(\frac{\alpha-1}{2\alpha-1})p>1$ since by hypothesis we have $p>\frac{3\alpha-2}{\alpha-1}$ and $1<\alpha<2$) and we finally obtain :
\begin{eqnarray*}
\|L(\vU)\|_{E} &=& \left\|\left(\begin{array}{c} \displaystyle{\int_0^t} \mathfrak{p}_{t-s}\ast\mathbb{P}(\theta \vec{e}_d)ds\\  0\end{array}\right)\right\|_E\leq \left\| \displaystyle{\int_0^t} \mathfrak{p}_{t-s}\ast\mathbb{P}(\theta \vec{e}_d)ds\right\|_{\mathcal{M}_\alpha^{p,\frac{d+\alpha}{\alpha-1}}}\leq  C \|\theta\|_{\mathcal{M}_\alpha^{\mathfrak{p},\frac{d+\alpha}{2\alpha-1}}}\\
&\leq & \frac{C}{\mathfrak{C}}\left(\|\vu\|_{\mathcal{M}_\alpha^{p,\frac{d+\alpha}{\alpha-1}}}+ \mathfrak{C}\|\theta\|_{\mathcal{M}_\alpha^{\mathfrak{p},\frac{d+\alpha}{2\alpha-1}}}\right)=\frac{C}{\mathfrak{C}}\|\vU\|_{E}.
\end{eqnarray*}
\item We will study now the bilinear term $B(\vU, \vU)$ and we have 
\begin{eqnarray}
\|B(\vU, \vU)\|_{E}&\leq &\underbrace{\left\| \displaystyle{\int_0^t} \mathfrak{p}_{t-s}\ast\mathbb{P}\left(\mathrm{div}(\vec{u}\otimes \vec{u})\right)ds\right\|_{\mathcal{M}^{p, \frac{d+\alpha}{\alpha-1}}_\alpha}}_{(1_B)}\notag\\
&&+\underbrace{\left\|  \displaystyle{\int_0^t} \mathfrak{p}_{t-s}\ast \mathrm{div}(\theta \vec{u})ds\right\|_{\mathcal{M}^{\mathfrak{p}, \frac{d+\alpha}{2\alpha-1}}_\alpha}}_{(2_B)}.\label{TermeBilineaire1}
\end{eqnarray}
For the term $(1_B)$ above we note that, since $\mathbb{P}\mathrm{div}(\cdot)$ is a vector of pseudo-differential operators with homogeneous symbol of degree $1$, smooth outside the origin, using the point-wise estimate of (\ref{Prop_estim_ptws}) for $\rho=1$, we obtain
\begin{eqnarray*}
\displaystyle{\int_0^t} \mathfrak{p}_{t-s}\ast\mathbb{P}\left(\mathrm{div}(\vec{u}\otimes \vec{u})\right)ds&\leq &C\displaystyle{\int_{\mathbb{R}\times \mathbb{R}^d}} \frac{1}{(|t-s|^\frac{1}{\alpha}+|x-y|)^{d+1}}|\vu(s,y)||\vu(s,y)|dyds\\[3mm]
&\leq & C\, \mathrm{I}_{\alpha-1}(|\vec{u}||\vec{u}|)(t,x),
\end{eqnarray*}
where we used the definition of the parabolic Riesz potentials given in (\ref{Def_ParabolicRiesz}) above. We then derive 
$$\left\| \displaystyle{\int_0^t} \mathfrak{p}_{t-s}\ast\mathbb{P}\left(\mathrm{div}(\vec{u}\otimes \vec{u})\right)ds\right\|_{\mathcal{M}^{p, \frac{d+\alpha}{\alpha-1}}_\alpha} \leq C \|\mathrm{I}_{\alpha-1}(|\vec{u}| |\vec{u}|)\|_{\mathcal{M}_\alpha^{p,\frac{d+\alpha}{\alpha-1}}}. $$
By the Proposition \ref{Prop_continuity_Riesz} with $\lambda=\frac{1}{2}$ we can write 
$$(1_B)\leq C\|\mathrm{I}_{\alpha-1}(|\vec{u}||\vec{u}|)\|_{\mathcal{M}_\alpha^{p,\frac{d+\alpha}{\alpha-1}}} \leq C \||\vec{u}||\vec{u}|\|_{\mathcal{M}_\alpha^{\frac{p}{2},\frac{d+\alpha}{2(\alpha-1)}}}\leq C  \|\vec{u}\|_{\mathcal{M}_\alpha^{p,\frac{d+\alpha}{\alpha-1}}} \|\vec{u}\|_{\mathcal{M}_\alpha^{p,\frac{d+\alpha}{\alpha-1}}},$$
where in the last estimate above we used the H\"older inequalities for parabolic spaces.\\
 
For the term $(2_B)$ of (\ref{TermeBilineaire1}), by the same pointwise arguments as above we have 
$$\displaystyle{\int_0^t} \mathfrak{p}_{t-s}\ast\mathbb{P}\left(\mathrm{div}(\theta\vec{u})\right)ds\leq C \,\mathrm{I}_{\alpha-1}(|\theta| |\vec{u}|),$$
from which we obtain (recall that $\mathfrak{p}=(\frac{\alpha-1}{2\alpha-1})p$):
$$\left\| \displaystyle{\int_0^t} \mathfrak{p}_{t-s}\ast\mathbb{P}\left(\mathrm{div}(\theta\vec{u})\right)ds\right\|_{\mathcal{M}^{\mathfrak{p}, \frac{d+\alpha}{2\alpha-1}}_\alpha} \leq C \|\mathrm{I}_{\alpha-1}(|\theta| |\vec{u}|)\|_{\mathcal{M}^{(\frac{\alpha-1}{2\alpha-1})p, \frac{d+\alpha}{2\alpha-1}}_\alpha}.$$
Now, applying Proposition \ref{Prop_continuity_Riesz} with $\lambda=\frac{2\alpha-1}{3\alpha-2}$ we obtain 
$$ \|\mathrm{I}_{\alpha-1}(|\theta| |\vec{u}|)\|_{\mathcal{M}^{\mathfrak{p}, \frac{d+\alpha}{2\alpha-1}}_\alpha}\leq C \||\theta| |\vec{u}|\|_{\mathcal{M}_\alpha^{(\frac{2\alpha-1}{3\alpha-2})\mathfrak{p},\frac{d+\alpha}{3\alpha-2}}}=C \||\theta| |\vec{u}|\|_{\mathcal{M}_\alpha^{(\frac{\alpha-1}{3\alpha-2})p,\frac{d+\alpha}{3\alpha-2}}},$$
since $\mathfrak{p}=(\frac{\alpha-1}{2\alpha-1})p$ we have $(\frac{2\alpha-1}{3\alpha-2})\mathfrak{p}=(\frac{\alpha-1}{3\alpha-2})p$. Note that we have the condition $1<(\frac{\alpha-1}{3\alpha-2})p$ (as by hypothesis we have $\frac{3\alpha-2}{\alpha-1}<p\leq \frac{d+\alpha}{\alpha-1}$) which is needed in order to apply the Proposition \ref{Prop_continuity_Riesz}.
\begin{Remarque}\label{RemarqueValueP}
Note here that the constraint $1<(\frac{\alpha-1}{3\alpha-2})p$  is crucial in order to perform our computations (\emph{i.e.} using the boundedness property of the parabolic Riesz transform in parabolic Morrey spaces via the Proposition \ref{Prop_continuity_Riesz}) and this forces the condition $p>\frac{3\alpha-2}{\alpha-1}>4$. The value of the parameter $p$ can not be made small using these arguments. 
\end{Remarque}

By the H\"older inequalities in the parabolic spaces with $\frac{3\alpha-2}{(\alpha-1)p}=\frac{1}{p}+\frac{2\alpha-1}{(\alpha-1)p}$ and $\frac{3\alpha-2}{d+\alpha}=\frac{\alpha-1}{d+\alpha}+\frac{2\alpha-1}{d+\alpha}$, we can write
$$(2_B)\leq C\|\mathrm{I}_{\alpha-1}(|\theta| |\vec{u}|)\|_{\mathcal{M}^{(\frac{\alpha-1}{2\alpha-1})p, \frac{d+\alpha}{2\alpha-1}}_\alpha}\leq C \|\vec{u}\|_{\mathcal{M}_\alpha^{p,\frac{d+\alpha}{\alpha-1}}} \|\theta\|_{\mathcal{M}_\alpha^{(\frac{\alpha-1}{2\alpha-1})p,\frac{d+\alpha}{2\alpha-1}}}.$$
With these two estimates for the terms $(1_B)$ and $(2_B)$, we can come back to the expression (\ref{TermeBilineaire1}) above to write (recall that $\mathfrak{p}=(\frac{\alpha-1}{2\alpha-1})p$)
\begin{eqnarray*}
\|B(\vU, \vU)\|_{E}&\leq & C  \|\vec{u}\|_{\mathcal{M}_\alpha^{p,\frac{d+\alpha}{\alpha-1}}} \|\vec{u}\|_{\mathcal{M}_\alpha^{p,\frac{d+\alpha}{\alpha-1}}}+C \|\vec{u}\|_{\mathcal{M}_\alpha^{p,\frac{d+\alpha}{\alpha-1}}} \|\theta\|_{\mathcal{M}_\alpha^{\mathfrak{p},\frac{d+\alpha}{2\alpha-1}}}\\
&\leq & C\|\vU\|_E\|\vU\|_E,
\end{eqnarray*}
which is the wished estimate.

\item For the external force $\vec{F}$, in order to establish the estimate (\ref{Estimations_Generales_ForcesExterieures}), we write:
\begin{eqnarray*}
\|\vec{F}\|_{E} = \left\|\left(\begin{array}{c} \displaystyle{\int_0^t} \mathfrak{p}_{t-s}\ast\mathbb{P}(\vec{f})ds\\
 \displaystyle{\int_0^t} \mathfrak{p}_{t-s}\ast g\,ds\end{array}\right)\right\|_E&= &\underbrace{\left\| \displaystyle{\int_0^t} \mathfrak{p}_{t-s}\ast\mathbb{P}(\vec{f})ds\right\|_{\mathcal{M}^{p, \frac{d+\alpha}{\alpha-1}}}}_{(1_F)}+\mathfrak{C}\underbrace{\left\| \displaystyle{\int_0^t} \mathfrak{p}_{t-s}\ast g\,ds\right\|_{\mathcal{M}^{\mathfrak{p}, \frac{d+\alpha}{2\alpha-1}}}}_{(2_F)},
\end{eqnarray*}
and we will study these two terms separately. For the term $(1_F)$ we write
$$ \displaystyle{\int_0^t} \mathfrak{p}_{t-s}\ast\mathbb{P}(\vec{f})ds=\displaystyle{\int_0^t} \mathfrak{p}_{t-s}\ast\mathbb{P}((-\Delta)^{\frac{\gamma}{2}}(-\Delta)^{-\frac{\gamma}{2}}\vec{f})ds,$$
and by the properties of the Leray projector we have 
$$ \displaystyle{\int_0^t} \mathfrak{p}_{t-s}\ast\mathbb{P}(\vec{f})ds=\displaystyle{\int_0^t} \mathbb{P}((-\Delta)^{\frac{\gamma}{2}}\mathfrak{p}_{t-s}\ast(-\Delta)^{-\frac{\gamma}{2}}\vec{f})ds.$$
We note now that the operator $\mathbb{P}((-\Delta)^{\frac{\gamma}{2}}\mathfrak{p}_{t-s}\ast(\cdot) )$ has a convolution kernel $\mathcal{K}_{t-s}$ and from Proposition \ref{Prop_estim_ptws} with $\rho=\gamma$ we can obtain the following pointwise decay 
$$\left|\mathcal{K}_{t-s}(x)\right| \leq \frac{C}{\left(|t-s|^\frac{1}{\alpha}+|x|\right)^{d+\gamma}},$$
we can thus write
\begin{eqnarray*}
\left| \int_0^t \mathfrak{p}_{t-s}\ast\mathbb{P}(\vec{f})ds\right|&\leq &C\int_{-\infty}^{+\infty}\int_{\mathbb{R}^d}  \frac{1}{(|t-s|^\frac{1}{\alpha}+|x-y|)^{d+\gamma}}|(-\Delta)^{-\frac{\gamma}{2}}\vec{f}(s,y)| dyds\\[2mm]
&\leq & C\mathrm{I}_{\alpha-\gamma}(|(-\Delta)^{-\frac{\gamma}{2}}\vec{f}|)(t,x),
\end{eqnarray*}
where we used the expression of the parabolic Riesz potential $\mathrm{I}_{\alpha-\gamma}$ given in (\ref{Def_ParabolicRiesz}). With this estimate at hand we have the control
$$\left\| \displaystyle{\int_0^t} \mathfrak{p}_{t-s}\ast\mathbb{P}(\vec{f})ds\right\|_{\mathcal{M}_\alpha^{p,\frac{d+\alpha}{\alpha-1}}} \leq C \left\| \mathrm{I}_{\alpha-\gamma}(|(-\Delta)^{-\frac{\gamma}{2}}\vec{f}|)\right\|_{\mathcal{M}_\alpha^{p,\frac{d+\alpha}{\alpha-1}}},$$
thus, applying Proposition \ref{Prop_continuity_Riesz} with $\lambda=\frac{\alpha-1}{2\alpha-1-\gamma}$, since we have by hypothesis 
$$0<\alpha-\gamma<(d+\alpha)\frac{2\alpha-1-\gamma}{d+\alpha}\quad \mbox{ and } \quad\mathfrak{m}=(\frac{\alpha-1}{2\alpha-1-\gamma})p,$$ 
we obtain the inequality
$$\left\| \mathrm{I}_{\alpha-\gamma}(|(-\Delta)^{-\frac{\gamma}   {2}}\vec{f}|)\right\|_{\mathcal{M}_\alpha^{p,\frac{d+\alpha}{\alpha-1}}} \leq C \left\| (-\Delta)^{-\frac{\gamma}{2}}\vec{f}\right\|_{\mathcal{M}_\alpha^{ \mathfrak{m},\frac{d+\alpha}{2\alpha-1-\gamma}}},$$
from which we easily deduce the bound (as we have $\mathfrak{r}=\frac{d+\alpha}{2\alpha-1-\gamma}$):
$$\left\| \displaystyle{\int_0^t} \mathfrak{p}_{t-s}\ast\mathbb{P}(\vec{f})ds\right\|_{\mathcal{M}_\alpha^{p,\frac{d+\alpha}{\alpha-1}}}\leq C \| \vec{f} \|_{\dot{\mathcal{W}}_\alpha^{-\gamma,\mathfrak{m},\frac{d+\alpha}{2\alpha-1-\gamma}}}=C\| \vec{f} \|_{\dot{\mathcal{W}}_\alpha^{-\gamma,\mathfrak{m},\mathfrak{r}}}.$$
For the second term $(2_F)$, we have by the same arguments: 
$$\left\| \displaystyle{\int_0^t} \mathfrak{p}_{t-s}\ast g\,ds\right\|_{\mathcal{M}^{\mathfrak{p}, \frac{d+\alpha}{2\alpha-1}}}\leq C \left\| \mathrm{I}_{\alpha-\delta}(|(-\Delta)^{-\frac{\delta}{2}}g|)\right\|_{\mathcal{M}_\alpha^{\mathfrak{p}, \frac{d+\alpha}{2\alpha-1}}},$$
thus, applying Proposition \ref{Prop_continuity_Riesz} with $\lambda=\frac{2\alpha-1}{3\alpha-1-\delta}$ we obtain 
$$ \left\| \mathrm{I}_{\alpha-\delta}(|(-\Delta)^{-\frac{\delta}{2}}g|)\right\|_{\mathcal{M}_\alpha^{\mathfrak{p}, \frac{d+\alpha}{2\alpha-1}}}\leq C  \left\|(-\Delta)^{-\frac{\delta}{2}} g\right\|_{\mathcal{M}_\alpha^{\mathfrak{n}, \frac{d+\alpha}{3\alpha-1-\delta}}},$$
where $\mathfrak{n}=(\frac{\alpha-1}{3\alpha-1-\delta})p$, which lead us to the control (since $\mathfrak{s}=\frac{d+\alpha}{3\alpha-1-\delta}$):
$$\left\| \displaystyle{\int_0^t} \mathfrak{p}_{t-s}\ast g\,ds\right\|_{\mathcal{M}^{\mathfrak{p}, \frac{d+\alpha}{2\alpha-1}}}\leq C\|g\|_{\dot{\mathcal{W}}_\alpha^{-\delta,\mathfrak{n},\mathfrak{s}}}.$$
With the estimates obtained for the terms $(1_F)$ and $(2_F)$ we finally have
\begin{eqnarray*}
\|\vec{F}\|_{E} &= &\left\| \displaystyle{\int_0^t} \mathfrak{p}_{t-s}\ast\mathbb{P}(\vec{f})ds\right\|_{\mathcal{M}^{p, \frac{d+\alpha}{\alpha-1}}}+\mathfrak{C}\left\| \displaystyle{\int_0^t} \mathfrak{p}_{t-s}\ast g\,ds\right\|_{\mathcal{M}^{\mathfrak{p}, \frac{d+\alpha}{2\alpha-1}}}\\
&\leq & C\|\vec{f} \|_{\dot{\mathcal{W}}_\alpha^{-\gamma,\mathfrak{m},\mathfrak{r}}}+C\mathfrak{C}\|g\|_{\dot{\mathcal{W}}_\alpha^{-\delta,\mathfrak{n},\mathfrak{s}}}.
\end{eqnarray*}
\end{itemize}
We have established the estimates (\ref{Estimations_Generales_DonneesInitiales})-(\ref{Estimations_Generales_TermeBilineaire}) needed to apply a fixed-point argument, and as the constants that appear in the computation do not depend on the time variable, under a smallness assumption over the initial data $(\vu_0, \theta_0)$ and on the external forces $(\vf, g)$ we can obtain a global in time mild solution for the fractional Boussinesq system.  \hfill $\blacksquare$
\section{Proof of the Theorem \ref{Theorem_Equivalences}}\label{Secc_Proof_Theorem_Equivalences}
We start noting that, if $(\vu, \theta)$ is a mild solution of the integral problem (\ref{Formulation_Integrale}) such that 
$$\vu \in \mathcal{M}_\alpha^{p,\frac{d+\alpha}{\alpha-1}}([0, +\infty[\times \mathbb{R}^d)\quad  \mbox{and} \quad \theta \in \mathcal{M}_\alpha^{\mathfrak{p},\frac{d+\alpha}{2\alpha-1}}([0, +\infty[\times \mathbb{R}^d),$$
then we have 
$$\|\mathfrak{p}_t\ast \vu_0\|_{ \mathcal{M}_\alpha^{p,\frac{d+\alpha}{\alpha-1}}}<+\infty \quad \mbox{and}\quad \|\mathfrak{p}_t\ast \theta_0\|_{\mathcal{M}_\alpha^{\mathfrak{p},\frac{d+\alpha}{2\alpha-1}}}<+\infty.$$
Thus, in order to prove the Theorem \ref{Theorem_Equivalences}, we will establish the following equivalence of norms
\begin{equation}\label{EquivParaMorreyTriebelLizorkin}
\|\mathfrak{p}_t\ast f\|_{\mathcal{M}^{p,q}_\alpha}\simeq \|f\|_{T^{-\frac{\alpha}{p}}_M{p,{\bf q}}},
\end{equation}
for a generic function $f:\mathbb{R}^d\longrightarrow \mathbb{R}$, where $1< p\leq q$ {$, p>\frac{\alpha}{d+\alpha}q$} and ${\bf q}$ is such that $\frac{d}{\bf q}=\frac{d+\alpha}{q}-\frac{\alpha}{p}$. Note in particular that with this definition of the parameter ${\bf q}$ (since $p\leq q$), we also have $p\leq {\bf q}$. Remark also that the condition $\frac{d}{\bf q}=\frac{d+\alpha}{q}-\frac{\alpha}{p}$ corresponds to the one given over the parameters of the Theorem \ref{Theorem_Equivalences} for the initial data $\vu_0$ and $\theta_0$.\\

The proof of the estimate $\|\mathfrak{p}_t\ast f\|_{\mathcal{M}^{p,q}_\alpha}\leq  \|f\|_{T^{-\frac{\alpha}{p}}_M{p,{\bf q}}}$ was already done in the proof of the Theorem \ref{Theorem_Triebel_Morrey} with some particular parameters. For the sake of completeness we give now a general version of this inequality and for this we write:
$$\|\mathfrak{p}_t\ast f\|_{\mathcal{M}_\alpha^{p,q}}=\underset{r>0}{\mathrm{sup}} \ \underset{(t,x)\in [0,+\infty[ \times \mathbb{R}^d}{\mathrm{sup}} \ \frac{1}{r^{(d+\alpha)(\frac{1}{p}-\frac{1}{q})}} \left(\displaystyle{\iint_{\{|t-s|^\frac{1}{\alpha}+|x-y|<r\}}}|\mathfrak{p}_s\ast f(y)|^pdyds \right)^\frac{1}{p},$$
Since the integration domain $\{|t-s|^\frac{1}{\alpha}+|x-y|<r\}$ is included in the set $[0,+\infty[\times B(x,r)$ we find the upper bound
\begin{eqnarray*}
\|\mathfrak{p}_t\ast f\|_{\mathcal{M}_\alpha^{p,q}}&\leq &\underset{r>0}{\mathrm{sup}} \ \underset{x\in \mathbb{R}^d}{\mathrm{sup}} \ \frac{1}{r^{(d+\alpha)(\frac{1}{p}-\frac{1}{q})}} \left(\displaystyle{\int_{B(x,r)}\int_{0}^{+\infty}}|\mathfrak{p}_s\ast f(y)|^pds dy\right)^\frac{1}{p}\\
&\leq & \underset{r>0}{\mathrm{sup}} \ \underset{x\in \mathbb{R}^d}{\mathrm{sup}} \ \frac{1}{r^{(d+\alpha)(\frac{1}{p}-\frac{1}{q})}} \frac{r^{d(\frac{1}{p}-\frac{1}{\bf q})}}{r^{d(\frac{1}{p}-\frac{1}{ \bf q})}}\left\|\left(\int_{0}^{+\infty}|\mathfrak{p}_s\ast f(y)|^pds\right)^{\frac{1}{p}}\right\|_{L^p(B(x,r))}.
\end{eqnarray*}
Since by hypothesis we have the relationship $\frac{d}{\bf q}=\frac{d+\alpha}{q}-\frac{\alpha}{p}$, we thus have the identity $(d+\alpha)(\frac{1}{p}-\frac{1}{q})=d(\frac{1}{p}-\frac{1}{\bf q})$ and we obtain the estimate
\begin{equation}\label{Estimation1Equivalence}
\|\mathfrak{p}_t\ast f\|_{\mathcal{M}_\alpha^{p,q}}\leq  \underset{r>0}{\mathrm{sup}} \ \underset{x\in \mathbb{R}^d}{\mathrm{sup}} \  \frac{1}{r^{d(\frac{1}{p}-\frac{1}{\bf q})}}\left\|\left(\int_{0}^{+\infty}|\mathfrak{p}_s\ast f(y)|^pds\right)^{\frac{1}{p}}\right\|_{L^p(B(x,r))}=\|f\|_{T^{-\frac{\alpha}{p}}_{M^{p,{\bf q}}}}.\\[5mm]
\end{equation}
We study now the reverse estimate $\|f\|_{T^{-\frac{\alpha}{p}}_M{p,{\bf q}}}\leq C\|\mathfrak{p}_t\ast f\|_{\mathcal{M}^{p,q}_\alpha}$ and to this end we write 
$$\|f\|_{T^{-\frac{\alpha}{p}}_{M^{p,{\bf q}}}}= \underset{r>0}{\mathrm{sup}} \ \underset{x\in \mathbb{R}^d}{\mathrm{sup}} \  \frac{1}{r^{d(\frac{1}{p}-\frac{1}{\bf q})}}\left(\int_{B(x,r)}\left(\int_{0}^{+\infty}|\mathfrak{p}_s\ast f(y)|^pds\right)dy\right)^{\frac{1}{p}},$$
by the Fubini theorem and since we have the identity $\frac{d}{\bf q}=\frac{d+\alpha}{q}-\frac{\alpha}{p}$ we can write
$$\|f\|_{T^{-\frac{\alpha}{p}}_{M^{p,{\bf q}}}}= \underset{r>0}{\mathrm{sup}} \ \underset{x\in \mathbb{R}^d}{\mathrm{sup}} \  \frac{1}{r^{(d+\alpha)(\frac{1}{p}-\frac{1}{q})}}\left(\int_{0}^{+\infty}\int_{B(x,r)}|\mathfrak{p}_s\ast f(y)|^p dyds\right)^{\frac{1}{p}},$$
from which we derive the control 
\begin{eqnarray}
\|f\|_{T^{-\frac{\alpha}{p}}_{M^{p,{\bf q}}}}&\leq & \underbrace{\underset{r>0}{\mathrm{sup}} \ \underset{x\in \mathbb{R}^d}{\mathrm{sup}} \  \frac{1}{r^{(d+\alpha)(\frac{1}{p}-\frac{1}{q})}}\left(\int_{0}^{r^\alpha}\int_{B(x,r)}|\mathfrak{p}_s\ast f(y)|^p dyds\right)^{\frac{1}{p}}}_{(A)}\notag\\
&&+ \underbrace{\underset{r>0}{\mathrm{sup}} \ \underset{x\in \mathbb{R}^d}{\mathrm{sup}} \  \frac{1}{r^{(d+\alpha)(\frac{1}{p}-\frac{1}{q})}}\left(\int_{r^\alpha}^{+\infty}\int_{B(x,r)}|\mathfrak{p}_s\ast f(y)|^p dyds\right)^{\frac{1}{p}}}_{(B)},\label{Estimation_EquivParaTriebel}
\end{eqnarray}
We remark now that, by the definition of the functional $\|\cdot\|_{\mathcal{M}^{p,q}_\alpha}$ given in the formula (\ref{Definition_MorreyPara}), the term $(A)$ above can be easily estimated in the following manner: 
$$(A)=\underset{r>0}{\mathrm{sup}} \ \underset{x\in \mathbb{R}^d}{\mathrm{sup}} \  \frac{1}{r^{(d+\alpha)(\frac{1}{p}-\frac{1}{q})}}\left(\int_{0}^{r^\alpha}\int_{B(x,r)}|\mathfrak{p}_s\ast f(y)|^p dyds\right)^{\frac{1}{p}}\leq \|\mathfrak{p}_s\ast f\|_{\mathcal{M}^{p,q}_\alpha}.$$
It remains now to study the term $(B)$ given in the expression (\ref{Estimation_EquivParaTriebel}) above and for this we write
$$\frac{1}{r^{(d+\alpha)(\frac{1}{p}-\frac{1}{q})}}\left(\int_{r^\alpha}^{+\infty}\int_{B(x,r)}|\mathfrak{p}_s\ast f(y)|^p dyds\right)^{\frac{1}{p}}\leq \frac{1}{r^{(d+\alpha)(\frac{1}{p}-\frac{1}{q})}}\left(\int_{r^\alpha}^{+\infty}\|\mathfrak{p}_s\ast f\|_{L^\infty}^pv_dr^d ds\right)^{\frac{1}{p}},$$
where $v_d$ is the volume of the $d$-dimensional unit ball.  We thus obtain
\begin{eqnarray*}
\frac{1}{r^{(d+\alpha)(\frac{1}{p}-\frac{1}{q})}}\left(\int_{r^\alpha}^{+\infty}\int_{B(x,r)}|\mathfrak{p}_s\ast f(y)|^p dyds\right)^{\frac{1}{p}}\leq C\frac{r^{\frac{d}{p}}}{r^{(d+\alpha)(\frac{1}{p}-\frac{1}{q})}}\left(\int_{r^\alpha}^{+\infty}s^{-\frac{\beta}{\alpha}p}\left(s^{\frac{\beta}{\alpha}}\|\mathfrak{p}_s\ast f\|_{L^\infty}\right)^p ds\right)^{\frac{1}{p}}\\
\leq C \underset{s>0}{\sup} \; s^{\frac{\beta}{\alpha}}\|\mathfrak{p}_s\ast f\|_{L^\infty} \quad \frac{r^{\frac{d}{p}}}{r^{(d+\alpha)(\frac{1}{p}-\frac{1}{q})}}\left(\int_{r^\alpha}^{+\infty}s^{-\frac{\beta}{\alpha}p} ds\right)^{\frac{1}{p}}.
\end{eqnarray*}
Recalling that  $\|f\|_{\dot{B}^{-\beta, \infty}_\infty}=\underset{s>0}{\sup} \; s^{\frac{\beta}{\alpha}}\|\mathfrak{p}_s\ast f\|_{L^\infty}$, fixing $\beta=\frac{d+\alpha}{q}=\frac{\alpha}{p}+\frac{d}{\bf{q}}$ and after integrating with respect to the time variable we obtain 
\begin{eqnarray*}
\frac{1}{r^{(d+\alpha)(\frac{1}{p}-\frac{1}{q})}}\left(\int_{r^\alpha}^{+\infty}\int_{B(x,r)}|\mathfrak{p}_s\ast f(y)|^p dyds\right)^{\frac{1}{p}}&\leq &C\|f\|_{\dot{B}^{-\beta, \infty}_\infty} \frac{r^{\frac{d}{p}}}{r^{(d+\alpha)(\frac{1}{p}-\frac{1}{q})}}r^{\frac{\alpha}{p}(1-\frac{\beta}{\alpha}p)}\\
&\leq & C \|f\|_{\dot{B}^{-\beta, \infty}_\infty} \frac{r^{\frac{d}{p}+\frac{\alpha}{p}-\beta}}{r^{(d+\alpha)(\frac{1}{p}-\frac{1}{q})}}\\
&\leq &C \|f\|_{\dot{B}^{-\beta, \infty}_\infty},
\end{eqnarray*}
since we have $\beta=\frac{d+\alpha}{q}$. Recall now that for a function $f:\mathbb{R}^d\longrightarrow \mathbb{R}$, for $z\in \mathbb{R}^d$ and for a real parameter $\lambda>0$, we have the identity $\|\mathfrak{p}_s\ast f(\lambda \cdot+z)\|_{\mathcal{M}^{p,q}_\alpha}=\lambda^{-\frac{d+\alpha}{q}}\|\mathfrak{p}_s\ast f\|_{\mathcal{M}^{p,q}_\alpha}$, from which we deduce the homogeneity and translation invariance for the quantity $\|\mathfrak{p}_s\ast \cdot\|_{\mathcal{M}^{p,q}_\alpha}$ thus, due to the maximality of the homogeneous Besov spaces (see \cite{Meyer}) we can derive the control 
\begin{equation}\label{EstimationHomogeneiteBesovMorrey}
\|f\|_{\dot{B}^{-\beta, \infty}_\infty}\leq C\|\mathfrak{p}_s\ast f\|_{\mathcal{M}^{p,q}_\alpha},
\end{equation}
since the quantities $\|\cdot\|_{\dot{B}^{-\beta, \infty}_\infty}$ and $\|\mathfrak{p}_s\ast \cdot\|_{\mathcal{M}^{p,q}_\alpha}$ share the same homogeneity and translation invariance (see a proof of this estimate in the appendix \ref{AppendixB} below). We thus easily obtain  
$$(B)= \underset{r>0}{\mathrm{sup}} \ \underset{x\in \mathbb{R}^d}{\mathrm{sup}} \  \frac{1}{r^{(d+\alpha)(\frac{1}{p}-\frac{1}{q})}}\left(\int_{r^\alpha}^{+\infty}\int_{B(x,r)}|\mathfrak{p}_s\ast f(y)|^p dyds\right)^{\frac{1}{p}}\leq  C \|\mathfrak{p}_t\ast f\|_{\mathcal{M}^{p,q}_\alpha},$$
and coming back to the expression (\ref{Estimation_EquivParaTriebel}) we finally have
\begin{equation}\label{Estimation2Equivalence}
\|f\|_{T^{-\frac{\alpha}{p}}_{M^{p,{\bf q}}}}\leq C'\|\mathfrak{p}_s\ast f\|_{\mathcal{M}^{p,q}_\alpha}.
\end{equation}
As we have the controls (\ref{Estimation1Equivalence}) and (\ref{Estimation2Equivalence}) we deduce the wished equivalence of norms (\ref{EquivParaMorreyTriebelLizorkin}): the proof of the Theorem \ref{Theorem_Equivalences} is now finished. \hfill $\blacksquare$
\appendix 
\section{Proof of the Proposition \ref{Propo_Inclusion_MorreyParabolique_Besov}}\label{AppendixA}

Assume that  $p,q$ are indexes such that $1\leq p<\frac{\alpha}{\beta}$ and $q=\frac{d+\alpha}{\beta}$ and consider a function $f:\mathbb{R}^d\longrightarrow \mathbb{R}$ such that $f\in \dot{B}^{-\beta, \infty}_\infty(\mathbb{R}^d)$. Since by the formula (\ref{Definition_MorreyPara}) above we have 
$$\|\mathfrak{p}_t\ast f\|_{\mathcal{M}_\alpha^{p,q}}=\underset{r>0}{\mathrm{sup}} \ \underset{(t,x)\in [0,+\infty[ \times \mathbb{R}^d}{\mathrm{sup}} \ \frac{1}{r^{(d+\alpha)(\frac{1}{p}- \frac{1}{q})}} \left(\displaystyle{\iint_{\{|t-s|^\frac{1}{\alpha}+|x-y|<r\}}}|\mathfrak{p}_s\ast f|^{p}dyds \right)^\frac{1}{p},$$
then we can write 
\begin{eqnarray*}
\|\mathfrak{p}_t\ast f\|_{\mathcal{M}_\alpha^{p, q}}=\underset{r>0}{\mathrm{sup}} \ \underset{(t,x)\in [0,+\infty[ \times \mathbb{R}^d}{\mathrm{sup}} \ \frac{1}{r^{(d+\alpha)(\frac{1}{p}-\frac{1}{q})}} \left(\displaystyle{\iint_{\{|t-s|^\frac{1}{\alpha}+|x-y|<r\}}}|s|^{-p\frac{\beta}{\alpha}}|s|^{p\frac{\beta}{\alpha}}|\mathfrak{p}_s\ast f|^{p}dyds \right)^\frac{1}{p}\\
\leq  \underset{s>0}{\sup}\; s^{\frac{\beta}{\alpha}}\|\mathfrak{p}_s\ast f\|_{L^\infty} \times \underset{r>0}{\mathrm{sup}} \ \underset{(t,x)\in [0,+\infty[ \times \mathbb{R}^d}{\mathrm{sup}} \ \frac{1}{r^{(d+\alpha)(\frac{1}{p}-\frac{1}{q})}} \left(\displaystyle{\iint_{\{|t-s|^\frac{1}{\alpha}+|x-y|<r\}}}|s|^{-p\frac{\beta}{\alpha}}dyds \right)^\frac{1}{p}.
\end{eqnarray*}
Thus, using the definition of the norm for the Besov spaces $\dot{B}^{-\beta, \infty}_\infty(\mathbb{R}^d)$ considered here (\emph{i.e.} we have $\|f\|_{\dot{B}^{-\beta, \infty}_\infty}= \underset{s>0}{\sup}\; s^{\frac{\beta}{\alpha}}\|\mathfrak{p}_s\ast f\|_{L^\infty}$), we can write
\begin{eqnarray} \label{Estimation_pourInclusionMorreyChinois}
\|\mathfrak{p}_t\ast f\|_{\mathcal{M}_\alpha^{p, q}}\leq \|f\|_{\dot{B}^{-\beta, \infty}_\infty}\times \underset{r>0}{\mathrm{sup}} \ \underset{(t,x)\in [0,+\infty[ \times \mathbb{R}^d}{\mathrm{sup}} \ \frac{1}{r^{(d+\alpha)(\frac{1}{p}-\frac{1}{q})}} \left(\displaystyle{\iint_{\{|t-s|^\frac{1}{\alpha}+|x-y|<r\}}}|s|^{-p\frac{\beta}{\alpha}}dyds \right)^\frac{1}{p},
\end{eqnarray}
and now we need to study the integrals above. Since the set $\{|t-s|^\frac{1}{\alpha}+|x-y|<r\}$ is included in the set   $\{s>0: |t-s|<r^\alpha\} \times \{y\in \mathbb{R}^d: |x-y|<r\}$, we have
$$I=\iint_{\{|t-s|^\frac{1}{\alpha}+|x-y|<r\}}|s|^{-p\frac{\beta}{\alpha}}dyds\leq C r^{d}\int_{\{|t-s|<r^\alpha\}}s^{-p\frac{\beta}{\alpha}}ds=Cr^{d}\int_{t-r^\alpha}^{t+r^\alpha}|s|^{-p\frac{\beta}{\alpha}}ds.$$
We decompose our study of the last integral in two cases. First, if $0\leq t\leq2r^\alpha$, the domain of integration is then included in the interval $[-r^\alpha, 3r^\alpha]$ and then it comes
$$I=Cr^d\int_{-r^\alpha}^{3r^\alpha}|s|^{-p\frac{\beta}{\alpha}}ds\leq Cr^d\int_{-3r^\alpha}^{3 r^\alpha}|s|^{-p\frac{\beta}{\alpha}}ds=2Cr^d\int_{0}^{3 r^\alpha}s^{-p\frac{\beta}{\alpha}}ds,$$
since $1\leq p<\frac{\alpha}{\beta}$, the previous integral is finite and we have 
$$I\leq C r^d r^{\alpha(1-p\frac{\beta}{\alpha})}=C r^{d+\alpha-p\beta}.$$
We consider now the case when $t>2r^\alpha$. Since we have here $t+r^\alpha>t-r^\alpha>r^\alpha>0$ and $-p\frac{\beta}{\alpha}<0$, we can write
$$I\leq Cr^{d}\int_{t-r^\alpha}^{t+r^\alpha}|s|^{-p\frac{\beta}{\alpha}}ds\leq Cr^{d}\int_{t-r^\alpha}^{t+r^\alpha}(r^\alpha)^{-p\frac{\beta}{\alpha}}ds=Cr^{d+\alpha-p\beta}.$$
With these estimates for $I$ we come back to the formula (\ref{Estimation_pourInclusionMorreyChinois}) and we obtain 
\begin{eqnarray*}
\|\mathfrak{p}_t\ast f\|_{\mathcal{M}_\alpha^{p,q}}&\leq &\|f\|_{\dot{B}^{-\beta, \infty}_\infty}\times \underset{r>0}{\mathrm{sup}} \ \underset{(t,x)\in [0,+\infty[ \times \mathbb{R}^d}{\mathrm{sup}} \ \frac{1}{r^{(d+\alpha)(\frac{1}{p}-\frac{1}{q})}} \left(\displaystyle{\iint_{\{|t-s|^\frac{1}{\alpha}+|x-y|<r\}}}|s|^{-p\frac{\beta}{\alpha}}dyds \right)^\frac{1}{p}\\
&\leq & C \|f\|_{\dot{B}^{-\beta, \infty}_\infty}\times\underset{r>0}{\mathrm{sup}} \ \underset{(t,x)\in [0,+\infty[ \times \mathbb{R}^d}{\mathrm{sup}} \ \frac{1}{r^{(d+\alpha)(\frac{1}{p}-\frac{1}{q})}} r^{\frac{d+\alpha}{p}-\beta}.
\end{eqnarray*}
But since we have $\frac{d+\alpha}{\beta}=q$, \emph{i.e.} $\beta=\frac{d+\alpha}{q}$, we can write 
\begin{eqnarray*}
\|\mathfrak{p}_t\ast f\|_{\mathcal{M}_\alpha^{p,q}}&\leq & C \|f\|_{\dot{B}^{-\beta, \infty}_\infty}\times\underset{r>0}{\mathrm{sup}} \ \underset{(t,x)\in [0,+\infty[ \times \mathbb{R}^d}{\mathrm{sup}} \ \frac{1}{r^{(d+\alpha)(\frac{1}{p}-\frac{1}{q})}} r^{\frac{d+\alpha}{p}-\frac{d+\alpha}{q}}=C \|f\|_{\dot{B}^{-\beta, \infty}_\infty},
\end{eqnarray*}
which is the wished inequality and this ends the proof of the Proposition \ref{Propo_Inclusion_MorreyParabolique_Besov}. \hfill$\blacksquare$
\section{A functional inequality}\label{AppendixB}
We will prove, for the sake of completeness, the estimate (\ref{EstimationHomogeneiteBesovMorrey}) with the proposition below:
\begin{Proposition}\label{Proposition_inclusion_Besovmax_parabolic_Morrey}
Let $d\geq 2$ and consider $1<\alpha < 2$ be a fixed parameter. Let $p,q$ be two real parameters such that $1\leq p\leq q < +\infty$ and define $\beta=\frac{d+\alpha}{q}$. Suppose that $f \in \mathcal{S}' / \mathcal{P}(\mathbb{R}^d)$ verifies $(t,x) \mapsto \mathds{1}_{t>0}\times\mathfrak{p}_t\ast f(x) \in \mathcal{M}_\alpha^{p,q}$. Then $f \in \dot{B}^{-\beta, \infty}_{\infty}(\mathbb{R}^d)$, and we have the estimate
$$\|f\|_{\dot{B}^{-\beta, \infty}_{\infty}} \leq C\|\mathfrak{p}_t\ast f\|_{\mathcal{M}_\alpha^{p,q}}.$$
\end{Proposition}
{\bf Proof.} It is enough to prove the case where $p=1$, since we have the continuous embedding $\mathcal{M}_\alpha^{p,q} \hookrightarrow \mathcal{M}_\alpha^{1,q}$. In order to simplify as much as possible the calculations, we will use the following norm on the space $\mathcal{M}_\alpha^{p,q}$, which is equivalent to the standard norm on this space. Namely, up to changing the constant $C$ in the desired estimate, we will use the norm
\begin{equation}\label{NormeMorreyParaboliqueEquiv}
\|\psi\|_{\mathcal{M}_\alpha^{p,q}} = \underset{r>0}{\sup} \ \underset{(t,x) \in \mathbb{R}\times \mathbb{R}^d}{\sup}\; \; \frac{1}{r^{(d+\alpha)(\frac{1}{p}-\frac{1}{q})}} \ \left( \displaystyle{\int_{\{|t-s|\leq r^\alpha\}}\int_{\{|x-y|_\infty \leq r\}}} |\psi(s,y)|^p dy ds \right)^\frac{1}{p},
\end{equation}
where $|z|_\infty=|(z_i)_{1\leq i \leq d}|_\infty = \underset{1\leq i \leq d}{\max} \ |z_i|$.\\

Now, consider $t>0$ and we express $\mathfrak{p}_t*f(x)$ as an integral in time by writing
$$\mathfrak{p}_t*f(x)= \frac{4}{t} \int_\frac{t}{4}^\frac{t}{2} \mathfrak{p}_t*f(x) ds.$$
Due to the semi-group property of the fractional heat kernel, for any $\frac{t}{4}\leq s \leq \frac{t}{2}$, we can write $\mathfrak{p}_t*f(x) = \mathfrak{p}_{t-s}*\left( \mathfrak{p}_s*f\right)(x)$ and we obtain
$$\mathfrak{p}_t*f(x) =  \frac{4}{t} \displaystyle{\int_\frac{t}{4}^\frac{t}{2}} \mathfrak{p}_{t-s}*\left(\mathfrak{p}_s*f\right)(x) ds,$$
from which we deduce 
\begin{eqnarray}
|\mathfrak{p}_t*f(x)| &\leq &\frac{4}{t} \int_\frac{t}{4}^\frac{t}{2} |\mathfrak{p}_{t-s}|*\left|\mathfrak{p}_s*f\right|(x) ds= \frac{4}{t} \displaystyle{\int_\frac{t}{4}^\frac{t}{2} \int_{\mathbb{R}^d}} |\mathfrak{p}_{t-s}(y)| |(\mathfrak{p}_s*f)(x-y)| dy ds. \label{estm_lemme_CN_Besov_1}
\end{eqnarray}
We study now the function $\mathfrak{p}_{t-s}(y)$ above. Recall that from homogeneity, we have the point-wise expression $\mathfrak{p}_{t-s}(y)=\frac{1}{(t-s)^\frac{d}{\alpha}}\mathfrak{p}_1\left(\frac{y}{(t-s)^\frac{1}{\alpha}}\right)$, furthermore, since $1< \alpha < 2$, the kernel $\mathfrak{p}_1$ possesses a certain decay at infinity guaranteeing integrability, namely one can find a constant $C=C(d,\alpha)$ such that
$$|\mathfrak{p}_1(y)|\leq \frac{C}{(1+|y|_\infty)^{d+1}},$$
(see \cite[Theorem 7.3.1, p. 320]{kolokoltsov}). Using this point-wise decay and inserting it in the expression of $\mathfrak{p}_{t-s}(y)$ yields
$$|\mathfrak{p}_{t-s}(y)| \leq \frac{1}{(t-s)^\frac{d}{\alpha}}\frac{C}{\left(1+\frac{|y|_\infty}{(t-s)^\frac{1}{\alpha}}\right)^{d+1}}.$$
Finally, since $\frac{t}{4} \leq s \leq \frac{t}{2}$, one has $\frac{t}{2}\leq t-s \leq \frac{3t}{4}$ and we get the following control 
$$|\mathfrak{p}_{t-s}(y)| \leq  \left(\frac{2}{t}\right)^\frac{d}{\alpha} \frac{C}{\left(1+\left(\frac{4}{3}\right)^\frac{1}{\alpha}\frac{|y|_\infty}{t^\frac{1}{\alpha}}\right)^{d+1}} \leq  \frac{1}{t^\frac{d}{\alpha}} \frac{C'}{\left(1+\frac{|y|_\infty}{t^\frac{1}{\alpha}}\right)^{d+1}}.$$
We then insert this estimate in the expression (\ref{estm_lemme_CN_Besov_1}) in order to obtain
$$|\mathfrak{p}_t*f(x)| \leq \frac{C''}{t^{1+\frac{d}{\alpha}}} \displaystyle{\int_\frac{t}{4}^\frac{t}{2} \int_{\mathbb{R}^d}} \frac{1}{\left(1+\frac{|y|_\infty}{t^\frac{1}{\alpha}}\right)^{d+1}} |(\mathfrak{p}_s*f)(x-y)| dy ds.$$
We now wish to use the information that $\mathfrak{p}_s*f$ belongs to $\mathcal{M}_\alpha^{1,q}$. To do this, we consider the covering of $\mathbb{R}^d$ consisting of the balls 
$$B_k=B_\infty(t^\frac{1}{\alpha}k, \tfrac{1}{2}t^\frac{1}{\alpha})=\{y \in \mathbb{R}^d: \ |y-t^\frac{1}{\alpha}k|_\infty \leq \tfrac{1}{2}t^\frac{1}{\alpha}\} = \{y \in \mathbb{R}^d: \ |\tfrac{y}{t^\frac{1}{\alpha}}-k|_\infty \leq \tfrac{1}{2} \},$$ 
for all $k \in \mathbb{Z}^d$ and for $t>0$. We can thus write
\begin{eqnarray}
|\mathfrak{p}_t*f(x)| & \leq & \frac{C}{t^{1+\frac{d}{\alpha}}} \displaystyle{\sum_{k\in \mathbb{Z}^d}} I_k,\label{estm_lemme_CN_Besov_2}
\end{eqnarray}
where we have
\begin{equation}\label{estm_lemme_CN_Besov_3}
I_k  =  \int_\frac{t}{4}^\frac{t}{2} \int_{B_k} \frac{1}{\left(1+\frac{|y|_\infty}{t^\frac{1}{\alpha}}\right)^{d+1}} |(\mathfrak{p}_s*f)(x-y)| dy ds.
\end{equation}
Now, for a fixed $k \in \mathbb{Z}^d$ and $y \in B_k$, observe that if $k=0$, we trivially have
$$\frac{1}{\left(1+\frac{|y|_\infty}{t^\frac{1}{\alpha}}\right)^{d+1}} \leq 1=\frac{1}{\left(1+|k|_\infty\right)^{d+1}},$$
while for $k \neq 0$, one has $|k|_\infty \geq 1$, hence by the triangular inequality we obtain
$$|\frac{y}{t^\frac{1}{\alpha}}|_\infty \geq |k|_\infty-|\frac{y}{t^\frac{1}{\alpha}}-k|_\infty \geq |k|_\infty-\frac{1}{2} \geq \frac{|k|_\infty}{2},$$
and we have
$$\frac{1}{\left(1+\frac{|y|_\infty}{t^\frac{1}{\alpha}}\right)^{d+1}} \leq \frac{1}{\left(1+\frac{|k|_\infty}{2}\right)^{d+1}} \leq \frac{C}{\left(1+|k|_\infty\right)^{d+1}}.$$
Thus, in all cases, we find by inserting this estimate in the integral in space of the expression (\ref{estm_lemme_CN_Besov_3}) above:
$$I_k \leq \frac{C}{\left(1+|k|_\infty\right)^{d+1}} \displaystyle{\int_\frac{t}{4}^\frac{t}{2} \int_{B_k}} |(\mathfrak{p}_s*f)(x-y)| dy ds.$$
We then apply the change of variables $y \mapsto x-y$ to obtain
$$I_k  \leq \frac{C}{\left(1+|k|_\infty\right)^{d+1}} \displaystyle{\int_\frac{t}{4}^\frac{t}{2} \int_{\widetilde{B}_{k,x}}} |(\mathfrak{p}_s*f)(y)| dy ds,$$
where $\widetilde{B}_{k,x}$ is the ball $B_\infty(x-t^\frac{1}{\alpha}k, \frac{1}{2} t^\frac{1}{\alpha})$. Note now that the integration domain $[\frac{t}{4}, \frac{t}{2}]\times B_\infty(x-t^\frac{1}{\alpha}k, \frac{1}{2} t^\frac{1}{\alpha})$ is included in the set $[-t, t]\times B_\infty(x-t^\frac{1}{\alpha}k, t^\frac{1}{\alpha})$, so we trivially obtain the bound
$$I_k  \leq  \frac{C}{\left(1+|k|_\infty\right)^{d+1}} \int_{\{|s-0|\leq t\}} \int_{\{|y-(x-t^\frac{1}{\alpha}k)|_\infty \leq t^\frac{1}{\alpha}\}} |\mathds{1}_{s>0}(\mathfrak{p}_s*f)(y)| dy ds.$$
We now use the assumption that $\mathfrak{p}_s*f$ belongs to the space $\mathcal{M}_\alpha^{1,q}$ (see the formula (\ref{NormeMorreyParaboliqueEquiv}) above) to write 
\begin{eqnarray*}
I_k & \leq &  \frac{C \left(t^\frac{1}{\alpha}\right)^{(d+\alpha)(1-\frac{1}{q})}}{\left(1+|k|_\infty\right)^{d+1}} \times\frac{1}{\left(t^\frac{1}{\alpha}\right)^{(d+\alpha)(1-\frac{1}{q})}} \int_{\{|s-0|\leq t\}} \int_{\{|y-(x-t^\frac{1}{\alpha}k)|_\infty \leq t^\frac{1}{\alpha}\}} |\mathds{1}_{s>0}(\mathfrak{p}_s*f)(y)| dy ds\\
&\leq &\frac{C \left(t^\frac{1}{\alpha}\right)^{(d+\alpha)(1-\frac{1}{q})}}{\left(1+|k|_\infty\right)^{d+1}} \|\mathfrak{p}_s*f\|_{\mathcal{M}_\alpha^{1,q}}= \frac{C}{\left(1+|k|_\infty\right)^{d+1}} \times t^{1+\frac{d}{\alpha}} \times t^{-\frac{d+\alpha}{\alpha  q}}\times \|\mathfrak{p}_s*f\|_{\mathcal{M}_\alpha^{1,q}}.
\end{eqnarray*}
We use the previous estimate of $I_k$ in the expression (\ref{estm_lemme_CN_Besov_2}) in order to obain 
\begin{eqnarray*}
|\mathfrak{p}_t*f(x)| & \leq & \frac{C}{t^{1+\frac{d}{\alpha}}} \displaystyle{\sum_{k\in \mathbb{Z}^d}} \frac{1}{\left(1+|k|_\infty\right)^{d+1}}\times t^{1+\frac{d}{\alpha}} \times t^{-\frac{d+\alpha}{\alpha q}}\times  \|\mathfrak{p}_s*f\|_{\mathcal{M}_\alpha^{1,q}}\\
&\leq& C  t^{-\frac{d+\alpha}{\alpha q}}  \|\mathfrak{p}_s*f\|_{\mathcal{M}_\alpha^{1,q}}
\displaystyle{\sum_{k\in \mathbb{Z}^d}} \frac{1}{\left(1+|k|_\infty\right)^{d+1}}=C t^{-\frac{\beta}{\alpha}} \|\mathfrak{p}_s*f\|_{\mathcal{M}_\alpha^{1,q}} \times\displaystyle{\sum_{k\in \mathbb{Z}^d}} \frac{1}{\left(1+|k|_\infty\right)^{d+1}} ,
\end{eqnarray*}
since we have $\beta=\frac{d+\alpha}{q}$. Note that the sum over $k \in \mathbb{Z}^d$ is finite, thus we have the control
$$|\mathfrak{p}_t*f(x)| \leq  C t^{-\frac{\beta}{\alpha}} \|\mathfrak{p}_s*f\|_{\mathcal{M}_\alpha^{1,q}},$$
from which we easily deduce the estimate
$$\underset{t>0}{\sup} \ t^\frac{\beta}{\alpha}\|\mathfrak{p}_t*f\|_{L^\infty}  \leq  C \|\mathfrak{p}_s*f\|_{\mathcal{M}_\alpha^{1,q}},$$
which is precisely the desired estimate $\|f\|_{\dot{B}^{-\beta, \infty}_{\infty}} \leq  C \|\mathfrak{p}_s*f\|_{\mathcal{M}_\alpha^{1,q}}$ and the proof of the Proposition \ref{Proposition_inclusion_Besovmax_parabolic_Morrey} is now ended. \hfill $\blacksquare$\\

\noindent {\bf Statement.} All authors have contributed to the manuscript substantially and have agreed to the final submitted version. This work does not have any conflict of interest.

\end{document}